\documentclass[aps,preprint,showpacs,superscriptaddress]{revtex4-1}  
\usepackage{graphicx}
\usepackage{subcaption}
\usepackage{bm}   
\usepackage{color}
\usepackage{amssymb}
\usepackage{amsmath}
\usepackage{float}

\begin{document}
	
	\title[Constant Bias and Weak Second Periodic Forcing : Tools to Mitigate Extreme Events]{Constant Bias and Weak Second Periodic Forcing : Tools to Mitigate Extreme Events}
	
	\author{S. Sudharsan}
	\affiliation{Department of Nonlinear Dynamics, Bharathidasan University, Tiruchirappalli - 620 024, Tamilnadu, India}
	\author{A. Venkatesan}
	\affiliation{Research Department of Physics, Nehru Memorial College (Autonomous), Puthanampatti, Tiruchirappalli - 621 007,  Tamil Nadu, India.}
	\author{M. Senthilvelan}
	\email[Correspondence to: ]{velan@cnld.bdu.ac.in}
	\affiliation{Department of Nonlinear Dynamics, Bharathidasan University, Tiruchirappalli - 620 024, Tamilnadu, India}
	\vspace{10pt}
	
	\begin{abstract}
		We propose two potentially viable non-feedback methods, namely (i) constant bias and (ii) weak second periodic forcing as tools to mitigate extreme events. We demonstrate the effectiveness of these two tools in suppressing extreme events in two well-known nonlinear dynamical systems, namely (i) Li\'enard system and (ii) a non-polynomial mechanical system with velocity dependent potential. As far as the constant bias is concerned, in the Li\'enard system, the suppression occurs due to the decrease in large amplitude oscillations and in the non-polynomial mechanical system the suppression occurs due to the destruction of chaos into a periodic orbit. As far as the second periodic forcing is concerned, in both the examples, extreme events get suppressed due to the increase in large amplitude oscillations. We also demonstrate that by introducing a phase in the second periodic forcing one can decrease the probability of occurrence of extreme events even further in the non-polynomial system. To provide a support to the complete suppression of extreme events, we present the two parameter probability plot for all the cases. In addition to the above, we examine the feasibility of the aforementioned tools in a parametrically driven version of the non-polynomial mechanical system. Finally, we investigate how these two methods influence the multistability nature in the Li\'enard system. 
	\end{abstract}
	
	%
	%
	%
	%
	%
	\maketitle
\section{Introduction}
Extreme events are rare events known for its aftermath destructive nature. These events emerge in the form of rogue waves, floods, cyclones, tsunamies, tornadoes, earthquakes, droughts, epidemics, epileptic seizures in physical, biological and natural systems \cite{Krause2015,epi,Sapsis2020,Ozturk2019}. In order to safeguard ourselves and the environment from the damage caused by these events, one has to devise various mitigation measures by investigating the mathematical models behind the phenomena. Here control measures represent various methods that one can implement in dynamical systems so that extreme events can be suppressed completely \cite{Jentsch2005}. In the literature, studies have been made on determining the emergence and mechanisms of extreme events in systems such as Li\'enard \cite{leo,kaviya}, microelectro mechanical system, solid $CO_2$ laser \cite{suresh}, loss-modulated $CO_2$ laser \cite{co2}, El Ni\~no southern oscillation \cite{elnino}, parametrically driven non-polynomial system \cite{sudharsan}, single and coupled Ikeda map \cite{ikeda}, coupled FitzHugh Nadumo model \cite{feudal1}, coupled Hindmarsh Rose model \cite{hr}, coupled Josephson junctions \cite{dana} and parametrically coupled chaotic populations \cite{ricker}. Although several works have been made to determine the mechanisms, only few studies have been devoted towards the determination of control measures. To name a few, we cite the following: (i) time delayed feedback is used to mitigate extreme events in Li\'enard system \cite{vkc}, (ii) threshold activated coupling has been used in coupled Ikeda map \cite{ikeda}, (iii) network mobiling was allowed in complex networks\cite{Chen2014}, (iv) closed-loop adaptive control was used in turbulent flows \cite{Farazmand2019}, (v) localized perturbations have been incorporated in spatially extended systems \cite{Bialonski2015}, (vi) influence of noise was utilized in optically injected lasers \cite{ZamoraMunt2014} and (vii) external forcing has been adopted in the parametrically driven non-polynomial system \cite{sudharsan}.

In the above studies, control measures have been devised using feedback mechanisms such as coupling or time delayed feedback. But rarely attempts have been made using non-feedback methods for the mitigation of extreme events. On the other hand, in the literature, non-feedback methods such as second weak periodic perturbation or constant bias have been utilized in various systems in a different context. For example, they have been used to control chaos in the driven Chua's circuit, Bonhoeffer-Van der Pol oscillator \cite{mlsr}, biharmonically driven sinusoidal potential system \cite{biharmo}, food chain model \cite{food}, electronic circuit system, the FitzHugh-Nagumo equation \cite{sr2}, and the Duffing equation \cite{control}. Non-feedback methods are easily implementable in experiments when compared to feedback methods. In the non-feedback methods, changes are not being made in the system according to its position in phase space. Instead changes are incorporated through the control parameters of the system \cite{mlsr}. By using non-feedback methods one can produce changes in the original system in such a way that the underlying dynamics is not totally deformed. Speed, high flexibility, non-requirement of online monitoring and online processing are some of the foremost advantages to consider non-feedback techniques \cite{sr2}. 

\par In one of our very recent works \cite{sudharsan}, we have shown that the first external force can suppress extreme events in the parametrically driven non-polynomial system. In Ref.~\cite{Kaveh2020}, a Markov-model-based chaos control was used to mitigate extreme events without the control of chaos in the system. To the best of our knowledge no report has been made on the control of extreme events in isolated systems either using constant bias or weak second periodic forcing. For the first time, we propose and demonstrate that these two methods are viable for the mitigation of extreme events. In particular, we test the applicability of these two methods on Li\'enard and a damped and driven non-polynomial mechanical system. Our results reveal that under the influence of constant bias not only extreme events get suppressed but also the chaotic nature is suppressed. Whereas in the case of second periodic forcing eventhough the extreme events get suppressed the chaos prevails. The outcome of our study is an important result as far as the study of mitigation of extreme events is concerned. Depending on the requirement, we can use either one of the two methods. Our results are promising in the field of nonlinear circuits, electro-optical system and so on.

\par In addition to the above, we also investigate the effect of inclusion of phase in the second periodic forcing. Here we find that extreme events get further suppressed in the non-polynomial system. Further, we also study the utility of these two methods in a parametrically driven version of the non-polynomial mechanical system. We also examine the multistability nature of the Li\'enard system under the influence of constant bias and second weak periodic forcing. Our results show that while the constant bias changes the multistability nature of the system, the second periodic forcing supresses the multistability nature only to a considerable extent. The suppression of multistability by these two methods has been proved using bifurcation plots, maximal Lyapunov exponent and basin of attraction. 

We organize our work as follows. In Sec. \ref{cons}, we investigate the effect of constant bias on extreme events in a Li\'enard system and in damped and driven non-polynomial system. In Sec. \ref{df}, we study the effect of second periodic forcing on extreme events in both the aforementioned systems. Further, we also examine the effect of phase on extreme events. We also analyze the combined effect of second periodic forcing and constant bias on these two systems. In Sec.~\ref{paramdfcons}, we study the effect of these non-feedback methods in a closely related dynamical system, namely parametrically driven non-polynomial system. In Sec.~{\ref{sec:multi}}, we study the effect of constant bias and double periodic forcing on the multistability nature of Li\'enard system. In Sec. \ref{discussion}, we consolidate the major outcomes of our study. Finally in Sec.~\ref{conclusion}, we conclude our work. 

\section{Influence of constant bias and mitigation of extreme events} 
\label{cons}
In this section, to begin, we investigate the influence of constant bias on Li\'enard system and on a non-polynomial mechanical system. We iterate both the systems for a very large time of about $10^9$ iterations. We carry out our analysis on extreme events after leaving sufficient transients. In particular, in the case of Li\'enard system we leave $10^5$ transients and in the case of non-polynomial system we leave $10^6$ transients. 
\subsection{Li\'enard system} 
\label{2.1}

The Li\'enard system with a constant bias $A$ is given by
\begin{eqnarray}
\dot{x}=y, \quad\dot{y}= -\alpha xy + \gamma x - \beta x^3 + f_1~\mathrm{cos}(\omega_1t)+A.
\label{eqn-l-cons}
\end{eqnarray}

Here, $\alpha$, $\beta$ and $\gamma$ are respectively the nonlinear damping, nolinearity strength and the internal frequency of the autonomous system. Further, $f_1$ represents the amplitude of the first external forcing, $\omega_1$ is the corresponding frequency and $A$ is the strength of constant bias. In our analysis, we fix the parameters as $\alpha=0.45$, $\beta=0.5$, $\gamma=0.5$, $f_1=0.2$, $\omega_1=0.7315$, $\omega_2=1.0$ and vary the parameter $A$. All these parameter values are chosen as in Ref. \cite{leo} such that the system exhibits extreme events due to interior crisis.

\begin{figure}[!ht]
	\centering
	\includegraphics[width=1.0\textwidth]{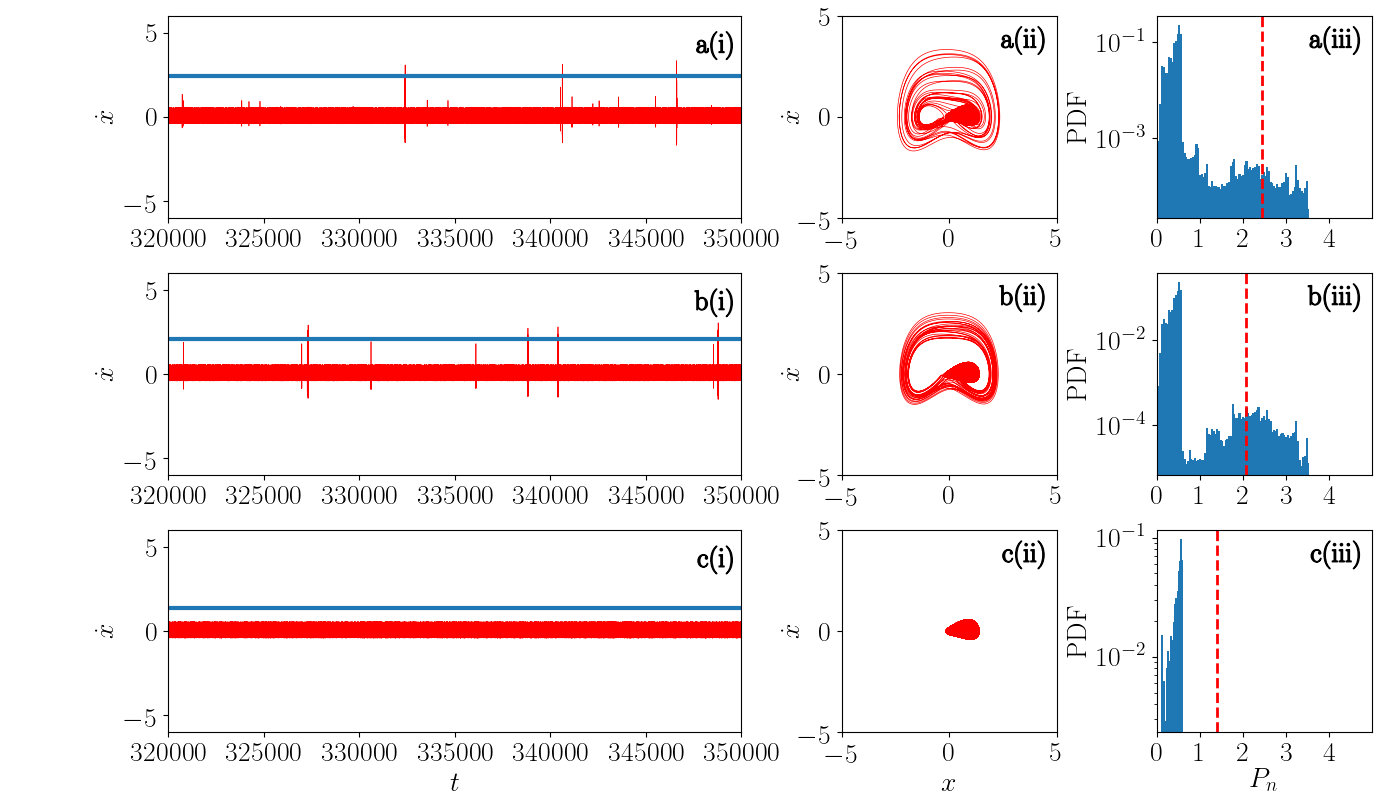}
	\caption{First panel represents the time series in $\dot{x}$, second panel represents the corresponding phase portrait and third panel represents the peak probability distribution function of system (\ref{eqn-l-cons}). Row (a) is for $A=0.0$, Row (b) is for $A=0.0001$ and Row (c) is for $A=0.001$. The value of other parameters are fixed as $\alpha=0.45$, $\beta=0.5$, $\gamma=0.5$, $f_1=0.2$, $\omega_1=0.7315$, $\omega_2=1.0$. }
	\label{l-ts-cons}
\end{figure}%

In dynamical systems, extreme events are determined using a predefined threshold value. This threshold value is called as qualifier threshold. The qualifier threshold is determined using the relation $x_{ee}=\langle x \rangle + n\sigma_x$ \cite{threshold}. Here $\langle x \rangle$ represents the average of all peaks and $\sigma_x$ corresponds to the standard deviation. The value of $n$ can be greater than or equal to four.  In the literature, in all the works on Li\'enard system in the context of extreme events, the value of $n$ in the qualifier threshold is fixed as eight. So, on par with the literature, we fix the value of $n$ as 8.

It is important to note that when $A=0$, under the influence of the first periodic forcing alone, system (\ref{eqn-l-cons}) exhibits extreme events at the point of interior crisis and at the point of intermittent expansion \cite{leo}. During the exhibition of extreme events, Li\'enard system's dynamics comprises of both large amplitude and small amplitude oscillations. The large amplitude oscillations are responsible for the extreme events \cite{leo,vkc}. When we introduce time-delayed feedback to the Li\'enard system and while increasing the feedback strength, extreme events were found to be controlled by the decrease in large amplitude oscillations \cite{vkc}. When there were no extreme events, the system displayed only small amplitude oscillations. 

Now we introduce a constant bias $A$, vary it and determine the effect of this constant bias on the extreme events. Figure~\ref{l-ts-cons} represents the time series in $\dot{x}$, corresponding phase portraits of the system and peak Probability Distribution Function (PDF) for various values of the constant bias. The peaks in column three of Fig.~\ref{l-ts-cons} are represented by $P_n$. The calculated threshold values are shown by a horizontal solid line in the time series plots and by a vertical dotted line in the PDF plots. Left panel represents the time series in $\dot{x}$, middle panel represents the phase portrait and the right panel represents the peak PDF. Figures \ref{l-ts-cons}a(i), a(ii) \& a(iii) show the dynamics during the occurrence of extreme events when $A=0$. In this state the dynamics exhibited by the system displays mixed mode oscillations comprising of both small and large amplitude oscillations. This fact can be seen both from the time series and phase portraits produced in Figs.~\ref{l-ts-cons}a(i) and a(ii) respectively. Further, the peak PDF shows a long tailed behaviour in such a way that peaks are present even beyond the threshold.  For a very small increase in the vlaue of $A$, say $A=0.0001$, we can observe that the number of large amplitude oscillations decreases considerably and as a consequence the number of extreme events also decreases. This can be verified from the Figs.~\ref{l-ts-cons}b(i), and  \ref{l-ts-cons}b(ii). However, the peak PDF still exhibits a long tailed behaviour (see Fig.~\ref{l-ts-cons}b(iii)).  Upon increasing the value of $A$ to $0.001$, the large amplitude oscillations are completely suppressed and only small amplitude oscillations prevail. This is visible from the time series shown in Fig.~\ref{l-ts-cons}c(i) where we can find not a single large amplitude oscillation. This fact also resembles in the phase portrait given in Fig.~\ref{l-ts-cons}c(ii). In this case the long tailed behaviour of the system vanishes and no peaks are seen beyond the threshold value, which is evident from Fig.~\ref{l-ts-cons}c(iii). In other words the extreme events are completely eliminated in the Li\'enard system. The large amplitude oscillations get suppressed upon increasing the strength of constant bias. This looks similar to that of supressing extreme events in the Li\'enard system with time-delayed feedback \cite{vkc}.

\begin{figure}[!ht]
	\centering
	\includegraphics[width=0.5\textwidth]{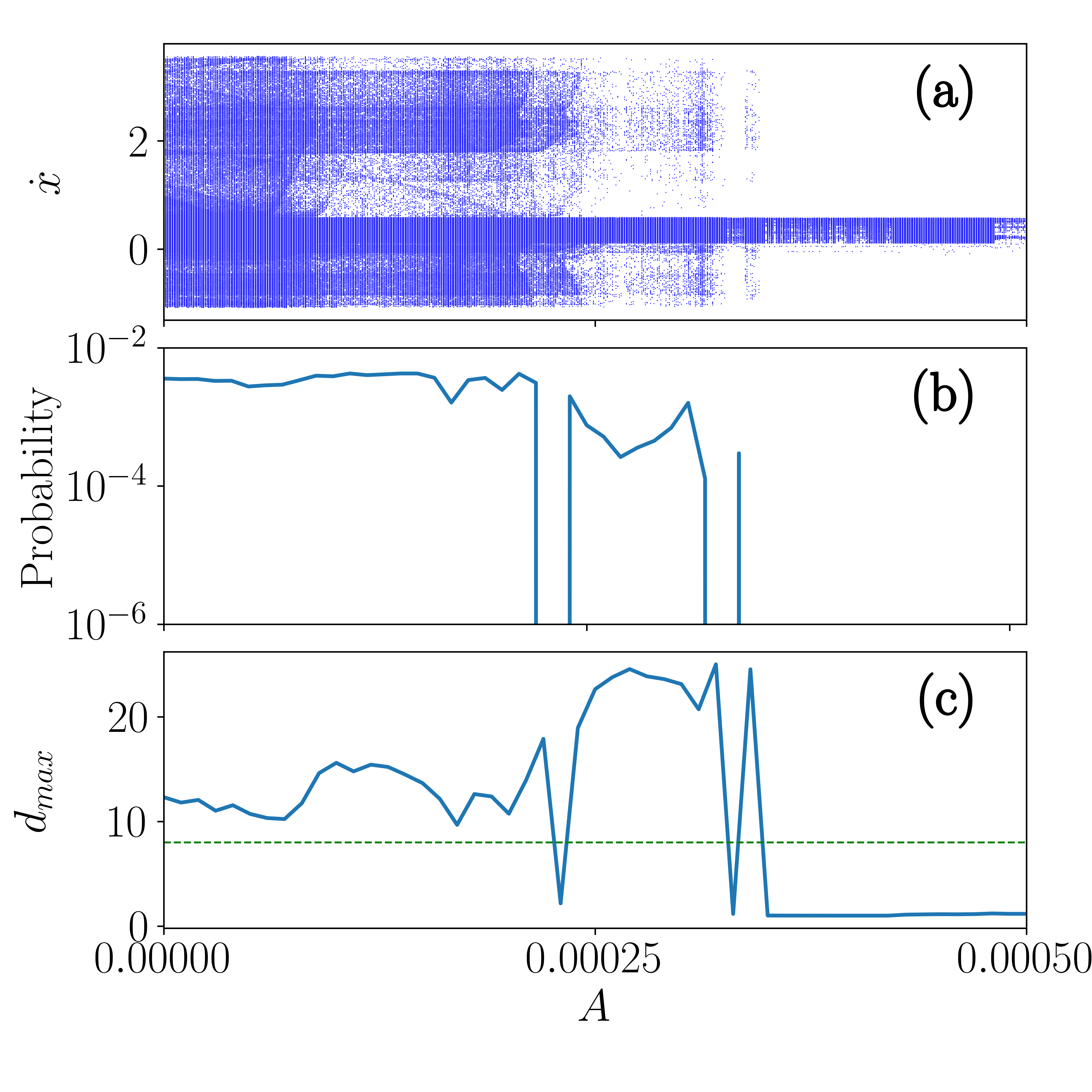}
	\caption{(a) represents Bifurcation diagram, (b) represents the probability of occurrence of extreme events and (c) represents the corresponding $d_{max}$ plot. The value of all the other parameters are as described in Fig.~\ref{l-ts-cons}.}
	\label{l-p-cons}
\end{figure}%

To further confirm the suppression of extreme events in the system (\ref{eqn-l-cons}), we plot the probability plot and the $d_{max}$ plot in Figs.~\ref{l-p-cons}(b) \& (c) respectively.  Here $d_{max}$ represents the number of standard deviations over the  average  value $\langle x\rangle$ corresponding to $\text{max}(x_{ee})$, where $\text{max}(x_{ee})$ represents the maximum value of the extreme event in the long time run.  The decreasing nature of extreme events with increasing constant bias is visible from Fig.~\ref{l-p-cons}(b). In particular, when $A=0.00034$, the probability of occurrence of extreme events reduces to zero as a result of transition of a large sized chaotic attractor to a small sized bounded chaos as shown in the bifurcation diagram in Fig.~\ref{l-p-cons}(a). That is, before the complete suppression of extreme events, the system displayed a large chaotic attractor comprising of both large and small amplitude oscillations. At the point where extreme events completely suppressed, only small bounded chaotic attractor prevails thereby zeroing the possibility of occurrence of extreme events. Here we have fixed the value of $n$ in the qualifier threshold as eight. The respective $d_{max}$ plot is shown in Fig.~\ref{l-p-cons}(c). It is calculated using the formula $d_{max}=\dfrac{\mathrm{max}(x_{ee})-\langle x\rangle}{\sigma_x}$ \cite{ikeda}. When we increase the value of $A$, the value of $d_{max}$  fluctuates and it exceeds the value eight whenever the extreme events occur. For regions where no extreme events occur, the value of $d_{max}$ lies below eight. The value of $d_{max}$ fluctuates in the extreme events occurring regimes and it becomes a constant in the non-extreme events regime. From the characterizations given in Fig.~\ref{l-p-cons}, we conclude that the constant bias can be used to control the extreme events in physical systems. 

\subsection{Non-polynomial system}
\label{npcons}

To further confirm the result regarding the suppression of extreme events by constant bias, we consider yet another dynamical system, a non-polynomial mechanical model which describes the motion of a freely sliding particle of unit mass on a parabolic wire rotating with a constant angular velocity $\omega_0$ about the axis $z=\sqrt{\lambda}x^2$. Here $1/\sqrt{\lambda}$ is the semi-latus rectum of the rotating parabola and $\lambda,~\omega_0>0$ \cite{nay}. With an  external forcing $f_1$ and a constant bias $A$, the equation of motion reads

\begin{eqnarray}
\dot{x}=y, \quad
\dot{y}=\frac{-\alpha y-\lambda x y^2 - \omega_0^2 x + f_1 \cos  \omega_e t + A}{1 + \lambda x^2}.
\label{eqn-m-cons}
\end{eqnarray} 

\begin{figure}[!ht]
	\centering
	\includegraphics[width=1.0\textwidth]{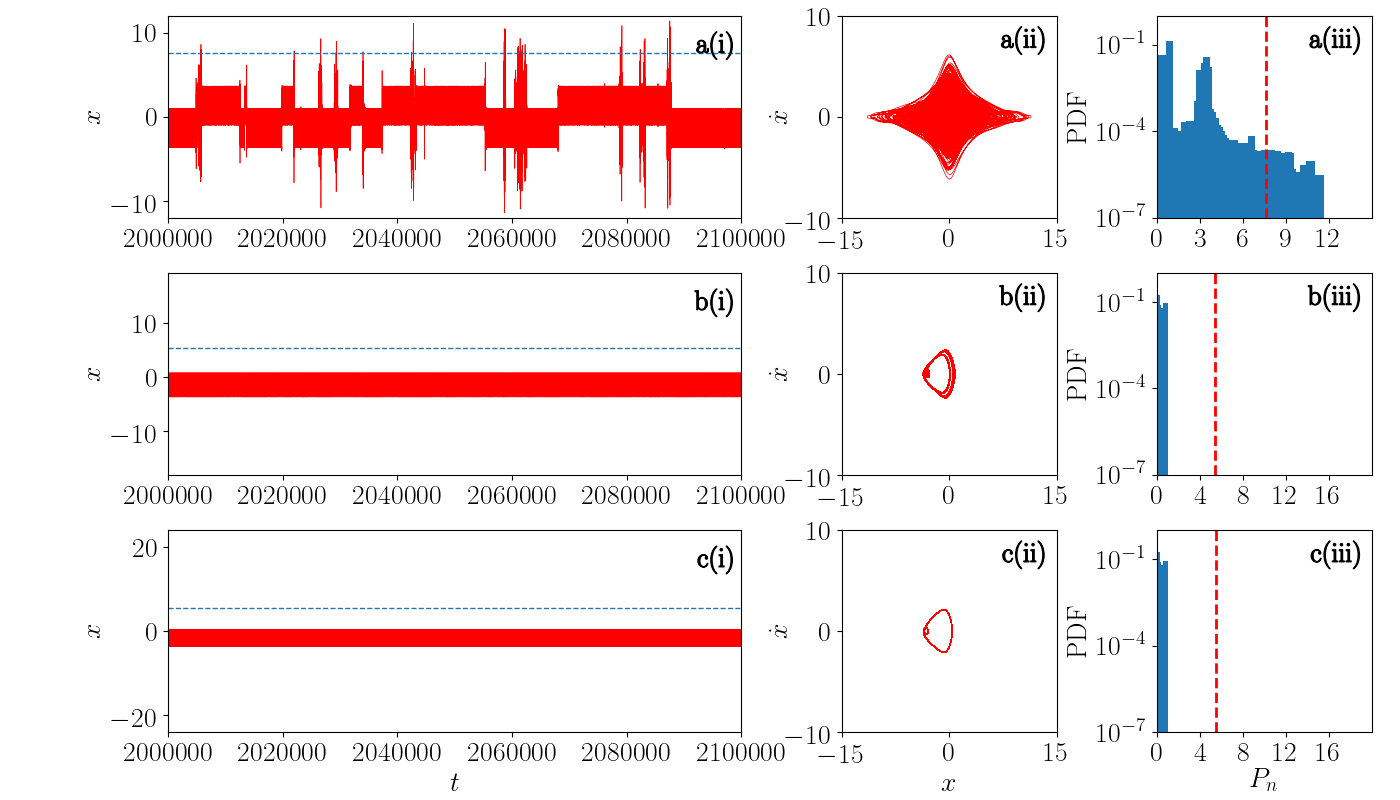}
	\caption{First panel represents the time series in $x$, second panel represents the corresponding phase portrait and third panel represents the peak probability distribution function of the non-polynomial mechanical systems under the influence of constant bias. Row (a) is for $A=0.0$, Row (b) is for $A=0.01$ and Row (c) is for $A=0.02$. The value of other parameters are  $\omega_0^2=0.25$, $\lambda=0.5$, $\alpha=0.2$, $f_1=3.1665$ and $\omega_1=1.0$.}
	\label{m-ts-cons}
\end{figure}

We fix the value of the parameters as $\omega_0^2=0.25$, $\lambda=0.5$, $\alpha=0.2$, $f_1=3.1665$ and $\omega_1=1.0$. These parameters are chosen such that system (\ref{eqn-m-cons}) exhibits extreme events with a single periodic forcing and $A=0$ \cite{sudharsan1}. Now we determine the influence of constant bias in suppressing extreme events in (\ref{eqn-m-cons}). We fix the value of $n$ in the qualifier threshold as four. The peaks are represented by $P_n$. In Fig.~\ref{m-ts-cons} we present the time series in $x$, phase portrait and  peak PDF of (\ref{eqn-m-cons}) in the first, second and third panels respectively. In this figure, first row is plotted for $A=0.0$, second row is plotted for $A=0.01$ and third row is plotted for $A=0.02$. When $A=0.0$, extreme events occurring in Eq.~(\ref{eqn-m-cons}) is confirmed from the trajectories that crosses the threshold (Fig.~\ref{m-ts-cons}a(i)), sparser trajectories at the edges of the phase portrait (Fig.~\ref{m-ts-cons}a(ii)) and the peak PDF which exhibits heavy tail distribution (Fig.~\ref{m-ts-cons}a(iii)). Upon adding a constant bias to the system and fixing its value be $A=0.01$, the number of extreme events occurring in the system completely reduces to zero. This can be confirmed from the time series plot given in Figs.~\ref{m-ts-cons}b(i) where now we see the absence of peaks beyond the threshold value. From Fig.~\ref{m-ts-cons}b(ii), we can infer that the chaotic nature of the attractor is now destroyed. Further, there is no long tail distribution in the peak PDF plot as we see in Fig.~\ref{m-ts-cons}b(iii). On further increasing the value of the constant bias to $A=0.02$, we observe that the system exhibits regular periodic behaviour. This is evident from the time series given in Fig.~\ref{m-ts-cons}c(i) and the corresponding phase portrait presented in Fig.~\ref{m-ts-cons}c(ii). Also no trajectory crosses the threshold value in the peak PDF which is evident from Fig.~\ref{m-ts-cons}c(iii). The suppression of extreme events occurs in (\ref{eqn-m-cons}) due to the destruction of the chaotic attractor and the birth of periodic orbit.

\begin{figure}[!ht]
	\centering
	\includegraphics[width=0.75\textwidth]{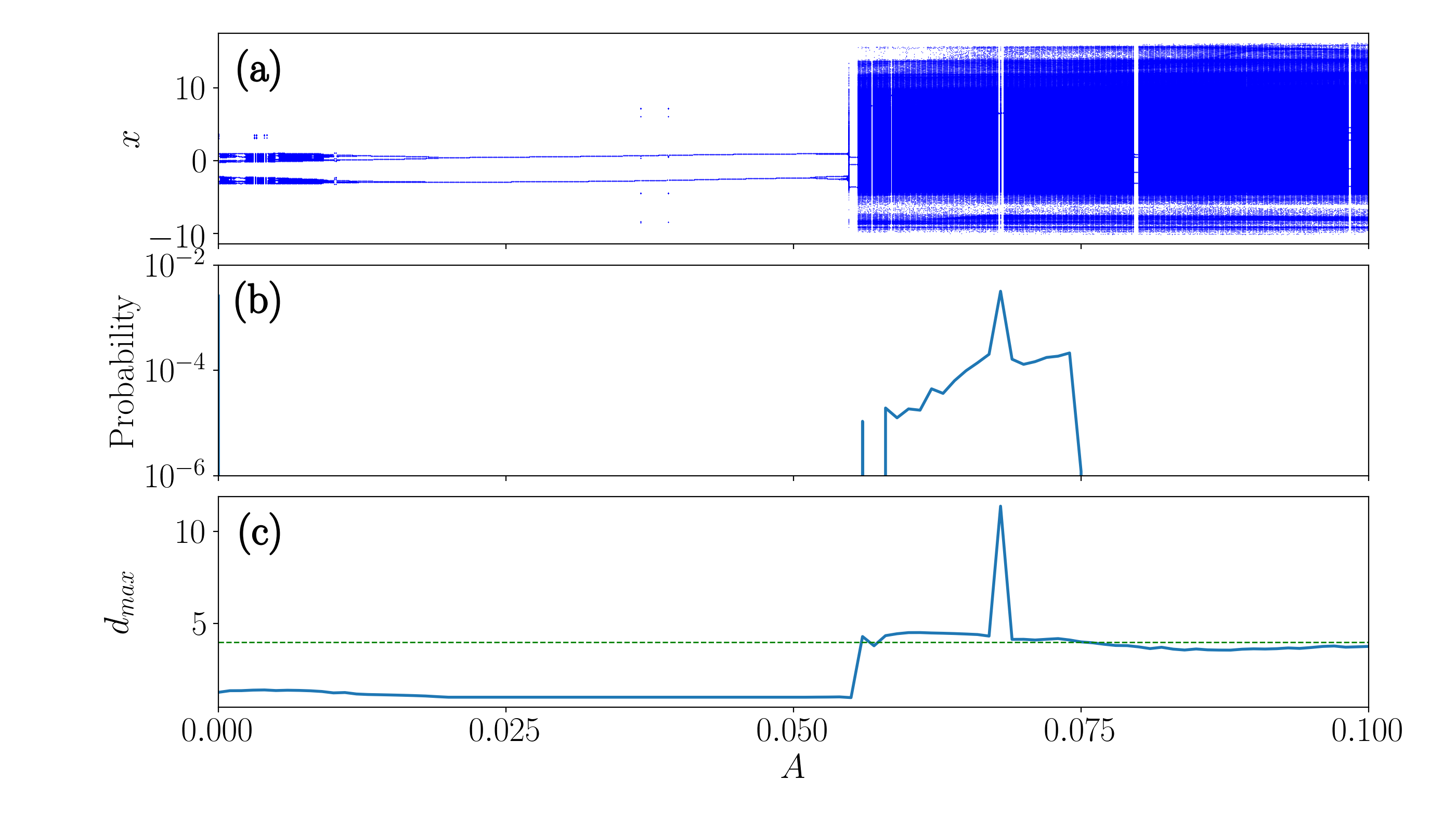}
	\caption{Plots of (a) Bifurcation diagram, (b) represents the probability of occurrence of extreme events and (c) represents the corresponding $d_{max}$. The values of all the other parameters are as mentioned in Fig.~\ref{m-ts-cons}.}
	\label{m-p-cons}
\end{figure}

In order to confirm the existence and suppression of extreme events, we plot the probability of occurrence of extreme events in Fig.~\ref{m-p-cons}. In Fig.~\ref{m-p-cons}(b), we can see the sudden reduction in the probability of occurrence of extreme events to zero. The chaotic attractor that was previously present takes the reverse period doubling route to exhibit a periodic behaviour. This can be confirmed from the bifurcation diagram shown in Fig~\ref{m-p-cons}(a). The periodic orbit sustains until $A=0.056$ where the chaotic attractor emerges again. After this, a short periodic window occurs. At $A=0.058$, the chaotic attractor remerges and the size of this chaotic attractor remains constant. At the point of re-emergence of this chaotic attractor, extreme events emerge and the probability of occurrence of extreme events increases gradually until $A=0.068$, after which the probability decreases and becomes zero at $A=0.075$. Thereafter extreme events get completely suppressed. The corresponding $d_{max}$ plot is shown in Fig.~\ref{m-p-cons}(c). Here the value of $d_{max}$ increases whenever extreme events occur and it crosses the value $n=4$. When no extreme events occur, $d_{max}$ lies below the value four. Although there exists a low probability for the occurrence of extreme events, when we increase the value of $A$, there is a considerable suppression. Thus, in the case of non-polynomial mechanical system under the influence of constant bias, extreme events are suppressed due to the destruction of chaotic attractor.  

\begin{figure}
	\centering
	\begin{subfigure}[b]{0.47\textwidth}
		\centering
		\includegraphics[width=\textwidth]{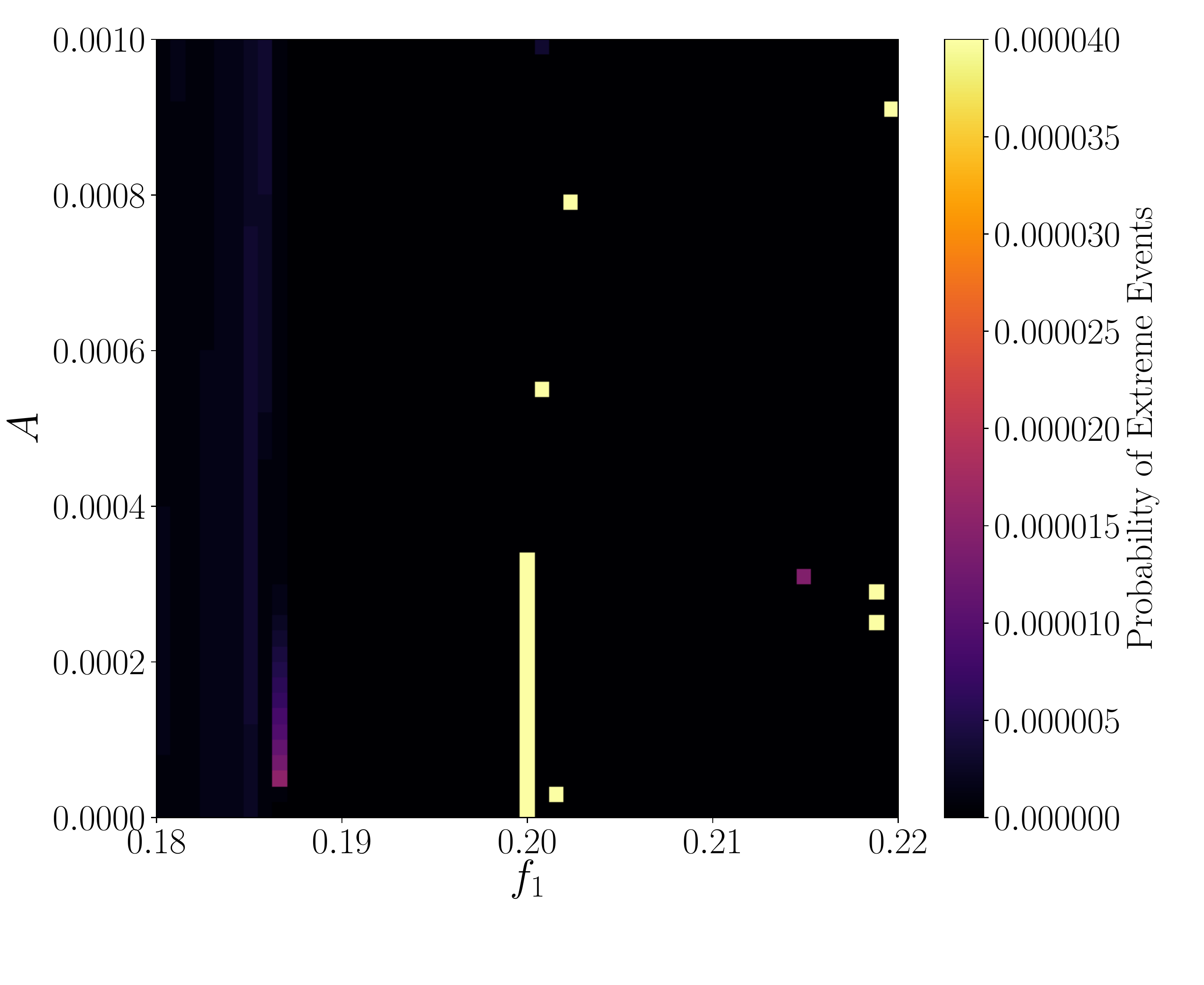}
		\caption{Li\'enard system}
		\label{two2a}
	\end{subfigure}
	\hfill
	\begin{subfigure}[b]{0.47\textwidth}
		\centering
		\includegraphics[width=\textwidth]{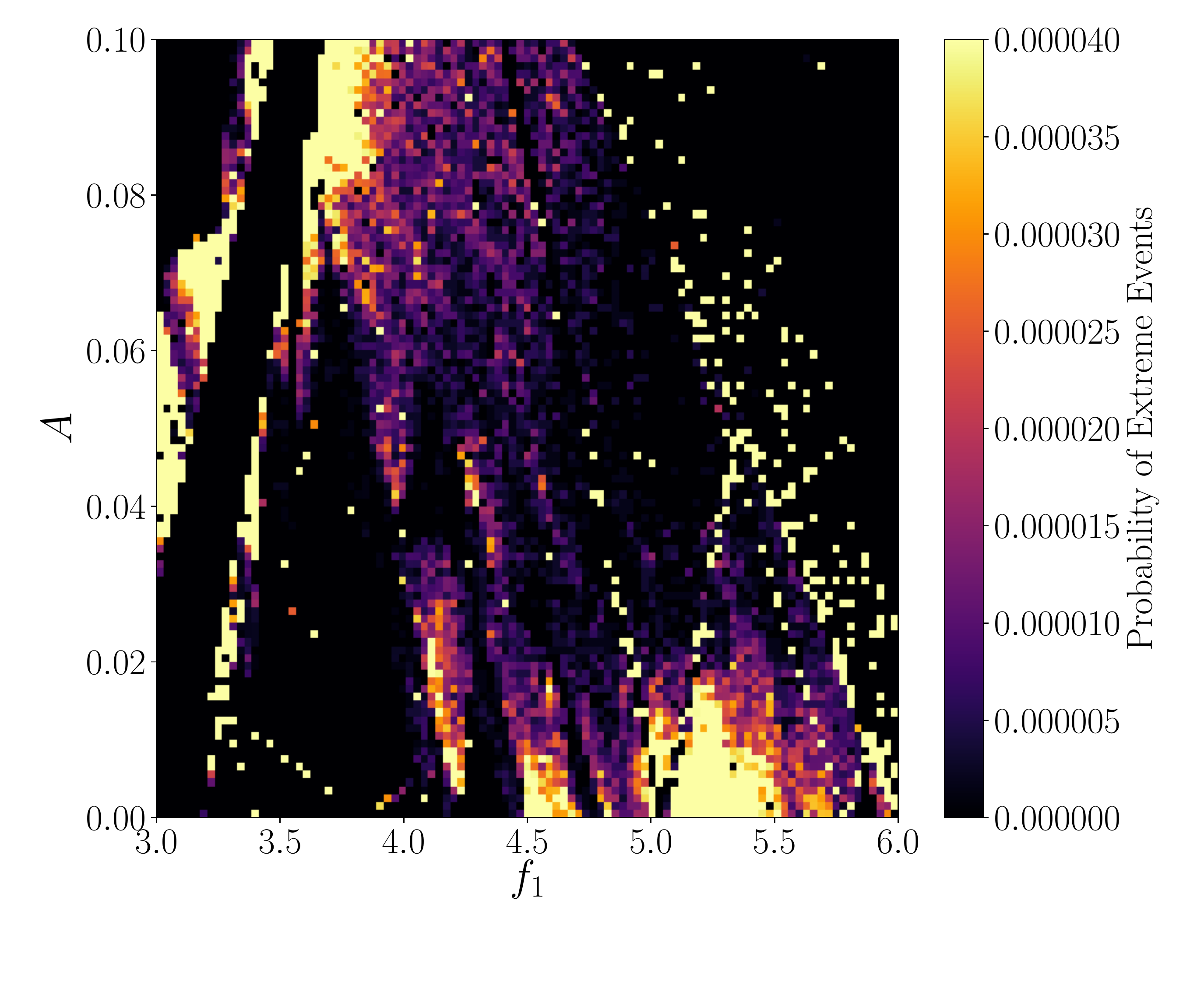}
		\caption{Mechanical system}
		\label{two2b}
	\end{subfigure}
	\caption{Probability of occurrence of extreme events in $f_1~-~f_2$ plane for the case of constant bias. The values of all the other parameters are as mentioned in Fig.~\ref{l-ts-cons} \& \ref{m-ts-cons}.}
	\label{two2}
\end{figure}

Finally, we present the two-parameter diagram in the $(f_1-A)$ plane for both systems (\ref{eqn-l-cons}) \& (\ref{eqn-m-cons}) with constant bias in Fig.~\ref{two2}. The colour panel in the right represents the value of the probability of occurrence of extreme events. In particular, black colour represents regions with zero extreme events and white colour represents regions with extreme event greater than the probability 0.00004, both in the case of Li\'enard system (Fig.~\ref{two2a}) and in the case of non-polynomial mechanical system (Fig.~\ref{two2b}). One can infer from the two parameter diagram that upon increasing the strength of constant bias the extreme events get suppressed.

\section{Influence of weak second periodic forcing and mitigation of extreme events}
\label{df}

In Sec.~\ref{cons}, we have shown that by augmenting constant bias one can suppress extreme events. In this section, we investigate the effect of weak second periodic forcing on extreme events.
\subsection{Li\'enard system}

The Li\'enard system with double periodic forcing is given by
\begin{eqnarray}
\dot{x}=y, \quad \dot{y}= -\alpha xy + \gamma x - \beta x^3 + f_1~\mathrm{cos}(\omega_1t)+f_2~\mathrm{cos}(\omega_2t).
\label{eqn-l-df}
\end{eqnarray}

The values of all the parameters in (\ref{eqn-l-df}) are same as given in the Sec.~\ref{2.1}. Instead of the constant bias  $A$, we introduce a second forcing of frequency $\omega_2$ with strength $f_2$. Now we vary it to determine how this weak second periodic forcing influences the extreme events.
\begin{figure}[!ht]
	\centering
	\includegraphics[width=1.0\textwidth]{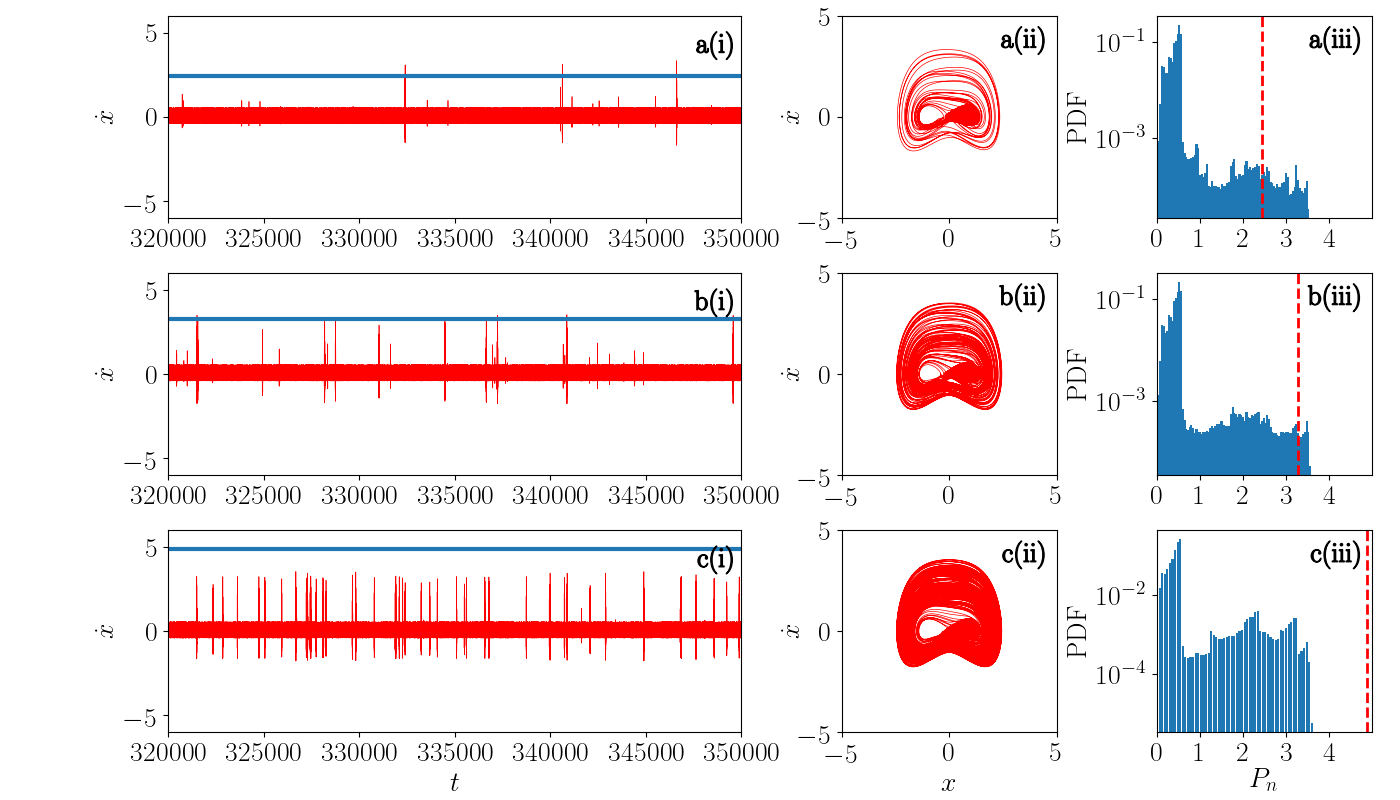}
	\caption{First panel represents the time series in $\dot{x}$, second panel represents the corresponding phase portrait and third panel represents the peak probability distribution function. Row (a) is for $f_2=0.0$, Row (b) is for $f_2=0.0001$ and Row (c) is for $f_2=0.001$. The values of all the other parameters are described in Fig.~\ref{l-ts-cons}.}
	\label{l-ts-df}
\end{figure}%

The time series in $\dot{x}$, phase portraits and the corresponing peak PDF for different values of the second forcing strength are shown in Fig.~\ref{l-ts-df}. For $f_2=0.0$, from Figs.~\ref{l-ts-df}a(i)\&a(ii), we can infer that the oscillations are in mixed mode comprising of both small and large amplitude oscillations. The occurrence of extreme events is evident from the trajectories that cross the threshold value (Fig.~\ref{l-ts-df}a(i)). In phase portrait, Fig.~\ref{l-ts-df}a(ii), the large amplitude oscillations contributing to the occurrence of extreme events can be clearly seen. In Fig.~\ref{l-ts-df}a(iii), we present the peak probability distribution.  The presence of long tail beyond the threshold value confirms the occurrence of extreme events. On increasing $f_2$ to $0.0001$, we observe that the number of large amplitude oscillations increases whereas the number of extreme events occurring in the given time domain reduces. This fact can be confirmed from the Figs.~\ref{l-ts-df}b(i) \& b(ii). As one visualizes, the number of extreme events get reduced even though there is an increase in the number of large amplitude oscillations. This is because of the increase in the number of large amplitude oscillations increases the threshold value thereby reducing the number of extreme events. The presence of extreme events is further confirmed by the presence of long tail from Fig.~\ref{l-ts-df}b(iii). But this time, in the peak PDF, the number of peaks that are present beyond the threshold value are comparatively lesser than the peaks that were present beyond the threshold in Fig.~\ref{l-ts-df}a(iii). This reflects the fact that the total number of extreme events occurring in the system decreases when we increase the second periodic forcing value.  On further increasing $f_2$ to $0.001$, the number of large amplitude oscillations increases thereby increasing the threshold value. This time the threshold value is well above the time series trajectories such that no trajectory crosses the threshold and hence ultimately no extreme events can be found. This can be confirmed from Fig.~\ref{l-ts-df}c(i). This is also evident from the dense phase portrait present in Fig.~\ref{l-ts-df}c(ii). The corresponding peak PDF also has no peaks beyond the threshold. Hence there is no long tail or heavy tail present in Fig.~\ref{l-ts-df}c(iii). This also confirms the absence of extreme events. All these facts substantiate that extreme events are suppressed in (\ref{eqn-l-df}) when we increase the value of the second periodic forcing strength $f_2$.

\begin{figure}[!ht]
	\centering
	\includegraphics[width=0.75\textwidth]{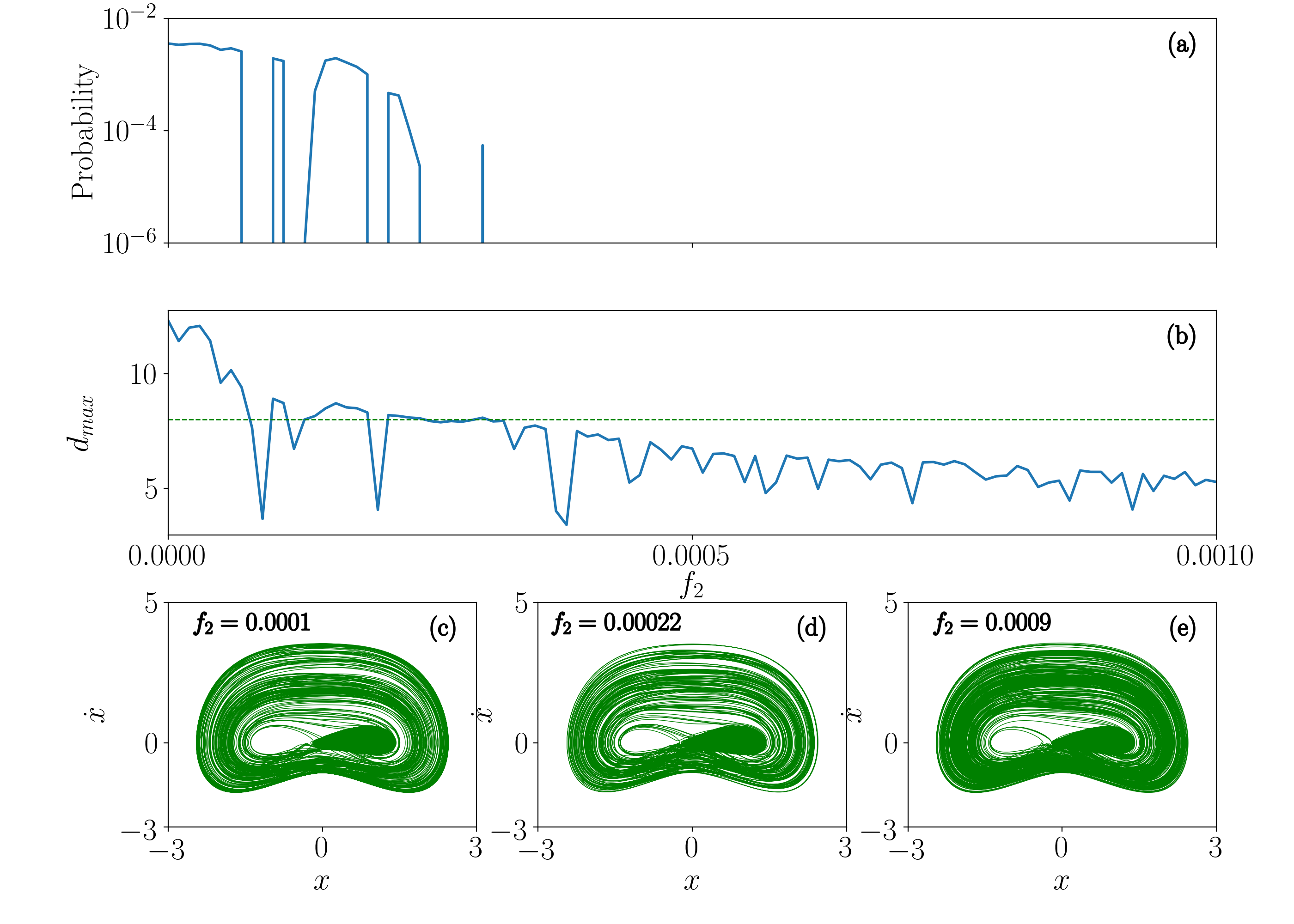}
	\caption{Plot (a) represents the Probability plot, (b) represents the $d_{max}$ plot. Plots (c), (d) \& (e) represents the phase portraits for three different forcing strengths $f_2=0.0001$, $f_2=0.00022$ and $f_2=0.0009$ respectively. The values of all the other parameters are described in Fig.~\ref{l-ts-cons}}
	\label{l-p-df}
\end{figure}%
In the above, we have demonstrated the suppression of extreme events by increasing $f_2$. Now, in order to substantiate our results, we plot, in Fig.~\ref{l-p-df}(a), the probability of occurrence of the extreme events and in Fig.~\ref{l-p-df}(b), the $d_{max}$ plot. We plot these two figures in order to confirm the suppression of extreme events. Phase portraits for three different values of the second periodic forcing strength are produced in Figs.~\ref{l-p-df}(c)-(e). This is to confirm whether the nature of trajectories are chaotic and also to verify whether large amplitude oscillations prevail even after the suppression of extreme events or not. The increasing nature of large amplitude oscillations are clearly exposed through the phase portraits. Figure~\ref{l-p-df}(a) shows the decrease in probability of occurrence of extreme events upon increasing the value of $f_2$. At a particular value of $f_2=0.0003$, the probability of occurrence of extreme events becomes zero, after which there is no signature of extreme events. Thus the second forcing $f_2$ suppresses the extreme events considerably even for a small value of $f_2$. Figure~\ref{l-p-df}(b) represents the $d_{max}$ plot. The value of $d_{max}$ decreases upon increasing $f_2$ and its value turns out to be greater than eight whenever the extreme events occur. In regions where no extreme events occur, the value of $d_{max}$ lies below eight. As the large amplitude oscillations increases, the value of $d_{max}$ decreases. The chaotic nature of the attractor is not destroyed upon the increasing value of $f_2$. From the outcome, we conclude that in the Li\'enard system one can suppress the extreme events by introducing a second periodic forcing and by increasing the strength of this forcing.

\subsection{Non-polynomial system}
\label{2.2}

Now we investigate the role of weak second periodic force on the non-polynomial mechanical model. The model with a weak second periodic force $f_2$ reads

\begin{eqnarray}
\dot{x}=y, \quad \dot{y}=\frac{-\alpha y-\lambda x y^2 - \omega_0^2 x + f_1 \cos  \omega_1 t +f_2 \cos  \omega_2 t}{1 + \lambda x^2}.
\label{eqn-m-df}
\end{eqnarray} 
The value of all the constants are same as given in Sec. \ref{npcons}, except $f_2=5.99865$. Let us introduce the second periodic forcing, as in the previous subsection, vary $f_2$ and study the problem at hand. In the present case, we fix the value of $n$ in the threshold qualifier as four. 

\begin{figure}[!ht]
	\centering
	\includegraphics[width=1.0\textwidth]{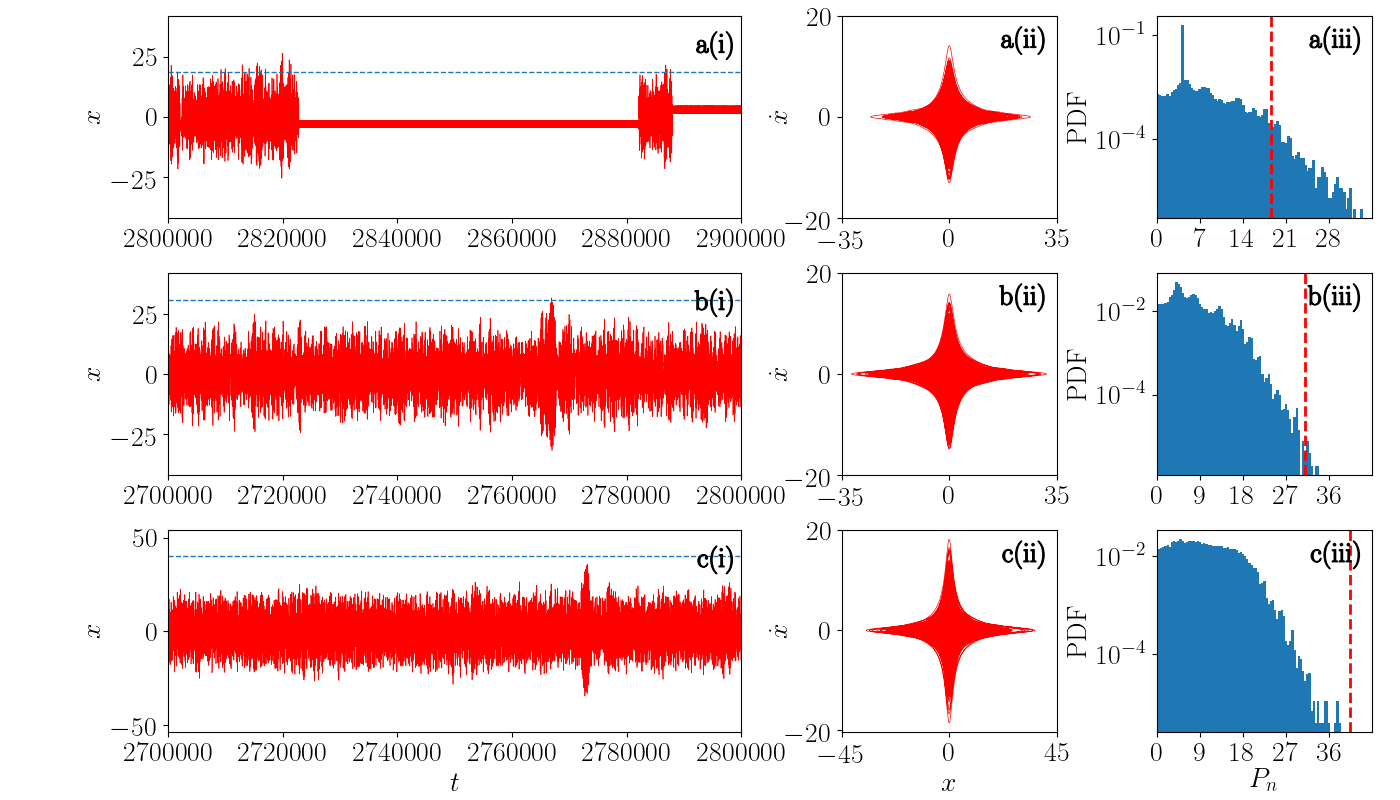}
	\caption{First panel represents the time series in $x$, second panel represents the corresponding phase portrait and third panel represents the peak probability distribution function. Row (a) is for $f_2=0.0$, Row (b) is for $f_2=0.04$ and Row (c) is for $f_2=0.8$. The values of all the other parameters are described in Fig.~\ref{m-ts-cons}.}
	\label{m-ts-df}
\end{figure}%

The occurrence of extreme events when $f_2=0.0$ is shown in Figs.~\ref{m-ts-df}a(i)-a(iii). It is evident from Fig.~\ref{m-ts-df}a(i) that more number of trajectories in the time series crosses the threshold value. In the correponding phase portrait given in Fig.~\ref{m-ts-df}a(ii), we can see sparser trajectories at the edges of the phase portrait which describe the extreme events. As far as the peak PDF is concerned, it exhibits long tailed behaviour. This is demonstrated in Fig.~\ref{m-ts-df}a(iii). When we increase the value of the second forcing strength $f_2$ to $0.04$, we can see the number of extreme events that occur in the system reduces considerably as the large amplitude oscillations of the chaotic attractor increases in the system. This is because the increase in large amplitude oscillations increases the value of threshold. This is clearly visible from the second panel in Fig.~\ref{m-ts-df}. Only one trajectory crosses the calculated threshold value as seen in the time series shown in Fig.~\ref{m-ts-df}b(i). The phase portrait in Fig.~\ref{m-ts-df}b(ii) has a lesser number of sparser trajectories when compared with Fig.~\ref{m-ts-df}a(ii). The peak PDF values also confirms this fact. On further increasing $f_2$ to 0.8, no extreme events occur and this is confirmed by the time series given in Fig.~\ref{m-ts-df}c(i) where we can see no trajectories crossing the threshold. Peak PDF in Fig.~\ref{m-ts-df}c(iii) also confirms the non-occurrence of extreme events as no peaks are found beyond the threshold.  
\begin{figure}[!ht]
	\centering
	\includegraphics[width=1.0\textwidth]{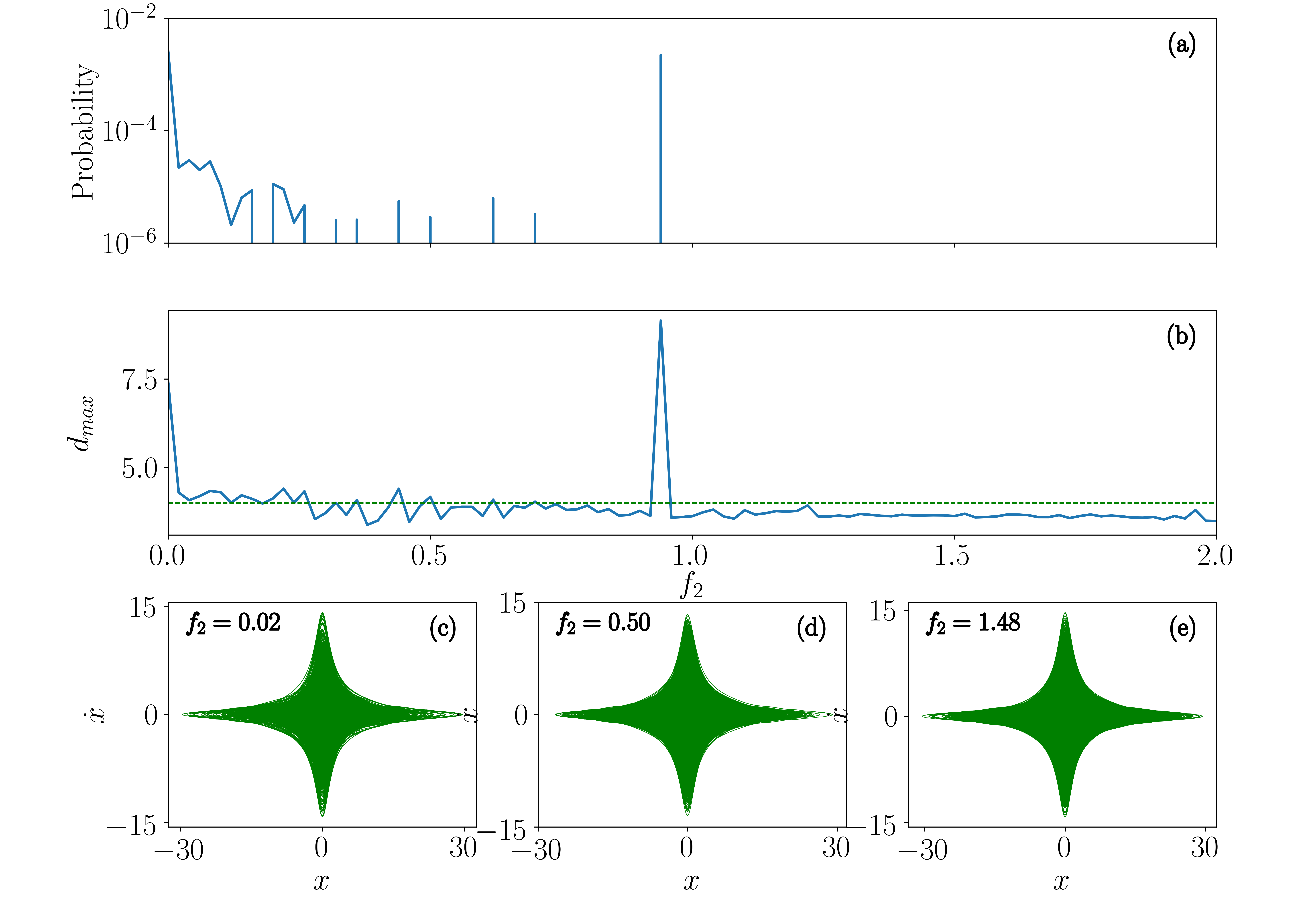}
	\caption{Plot (a) represents the Probability plot, (b) represents the $d_{max}$ plot. Plots (c), (d) \& (e) represents the phase portraits for three different forcing strengths $f_2=0.02$, $f_2=0.50$ and $f_2=1.42$ respectively. The values of all the other parameters are described in Fig.~\ref{m-ts-cons}.}
	\label{m-p-df}
\end{figure}%

\par In order to confirm the existence and suppression of extreme events under the influence of the second periodic forcing strength, we plot the probability of occurrence of extreme events in Fig.~\ref{m-p-df}(a). From Fig.~\ref{m-p-df}(a), we can infer that the probability of occurrence of extreme events decreases on increasing the second forcing strength. Beyond $f_2=0.94$, there is no occurrence of extreme events. This can also be confirmed from the $d_{max}$ plot which is given in Fig.~\ref{m-p-df}(b). Whenever extreme events occur, the value of $d_{max}$ rises above four and for non-extreme events regime it lies below four. We plot the phase portraits in Figs.~\ref{m-p-df}(c)-(e) for three different values of $f_2$. In these phase portraits one may notice that the nature of the attractor does not change. This confirms that the nature of the trajectory is chaotic and that it persists even after the suppression of extreme events. From the above outcome, we conclude that the second periodic forcing suppresses the occurrence of extreme events in the considered non-polynomial mechanical system also. 

\begin{figure}
	\centering
	\begin{subfigure}[b]{0.47\textwidth}
		\centering
		\includegraphics[width=\textwidth]{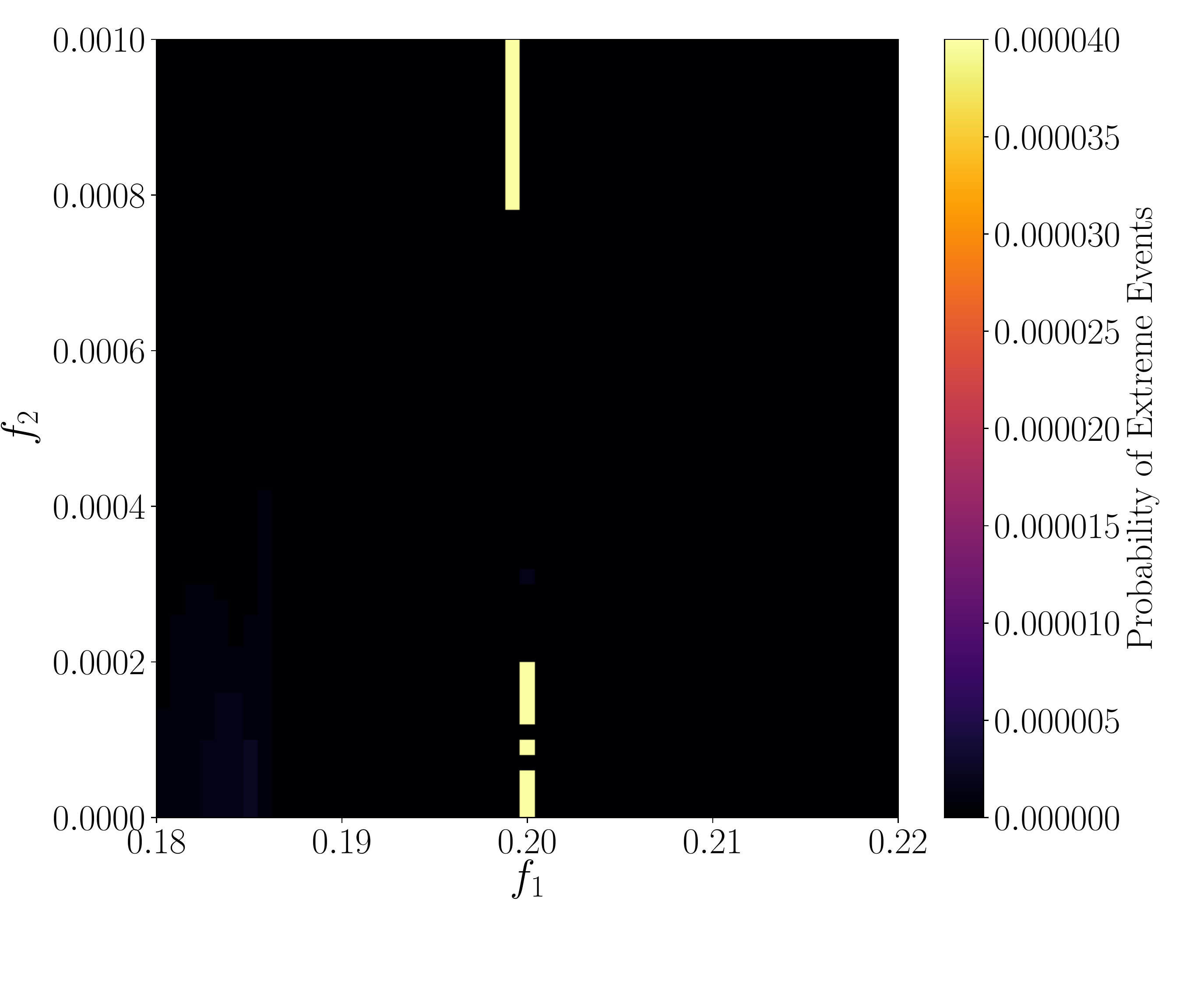}
		\caption{Li\'enard system}
		\label{two1a}
	\end{subfigure}
	\hfill
	\begin{subfigure}[b]{0.47\textwidth}
		\centering
		\includegraphics[width=\textwidth]{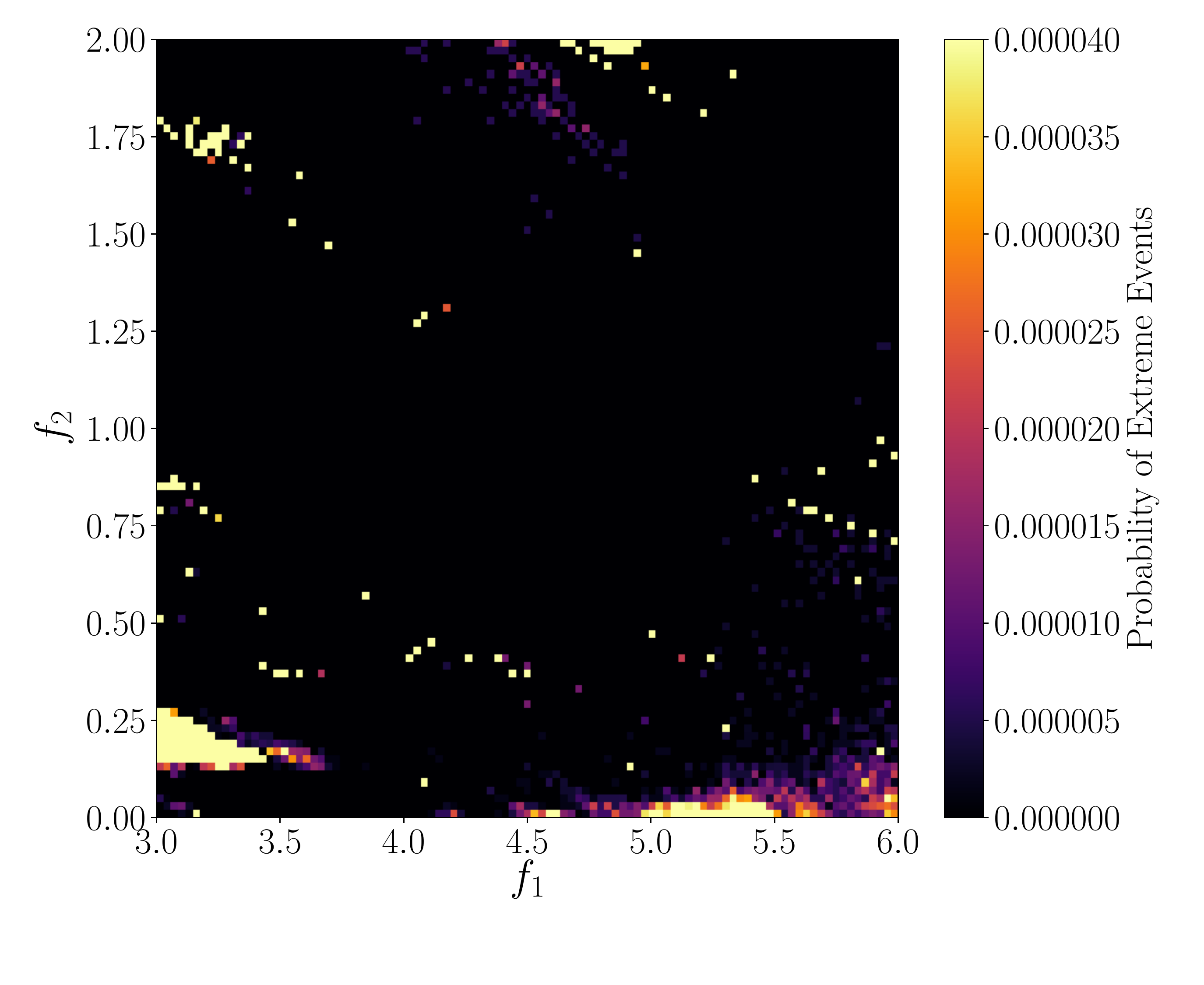}
		\caption{Mechanical system}
		\label{two1b}
	\end{subfigure}
	\caption{Probability of occurrence of extreme events in $f_1~-~f_2$ plane for the case of double forcing. The values of all the other parameters are described in Fig.~\ref{l-ts-cons} \& \ref{m-ts-cons}.}
	\label{two1}
\end{figure}

The two-parameter diagram in the $(f_1-f_2)$ plane for both the Li\'enard and non-polynomical mechanical system are given in Fig.~\ref{two1} in the presence of second periodic forcing. In this figure, the colour panel in the right represents various probabilities for the occurrence of extreme events. In particular, the apex colour, black represents region with zero extreme events and white represents region with extreme events greater than the probability 0.00004 for both the Li\'enard and the non-polynomial mechanical system. As one notices, most of the regions in Figs.~\ref{two1a} \& b are black and only for very few parameter values, extreme events exist. At maximum number of places the probability of occurrence of extreme events is zero. In particular, when the second forcing strength $f_2$ is increased, extreme events get mitigated completely.

As in the case of Li\'enard system there is no suppression of large amplitude oscillations in the non-polynomial mechanical system when it is under the influence of weak second forcing. But still the extreme events get suppressed. 

From the above discussion, we conclude that second periodic forcing also plays a crucial role in suppressing extreme events in physical systems. Since extreme events in most of the cases needs to be avoided, second periodic forcing can be considered as a viable tool to mitigate extreme events. 

\subsection{Influence of phase of the weak periodic force on the extreme events}
\label{dfcons1}
From Fig.~\ref{m-p-df}(a), we find that while increasing $f_2$, the extreme events emerge discontinuously. In order to check whether those extreme events can be suppressed or not, additionally, we introduce a phase $\phi$ in the second forcing term in the form of $\mathrm{cos}~(\omega_2~t + \phi)$ in Eq.~(\ref{eqn-m-df}). The underlying equation now reads

\begin{eqnarray}
\dot{x}=y,\quad \dot{y}=\frac{-\alpha y-\lambda x y^2 - \omega_0^2 x + f_1 \cos  \omega_1 t +f_2 \cos  (\omega_2 t+\phi)}{1 + \lambda x^2}.
\label{eqn-m-df-phi}
\end{eqnarray} 

\begin{figure}[!ht]
	\centering
	\includegraphics[width=0.75\textwidth]{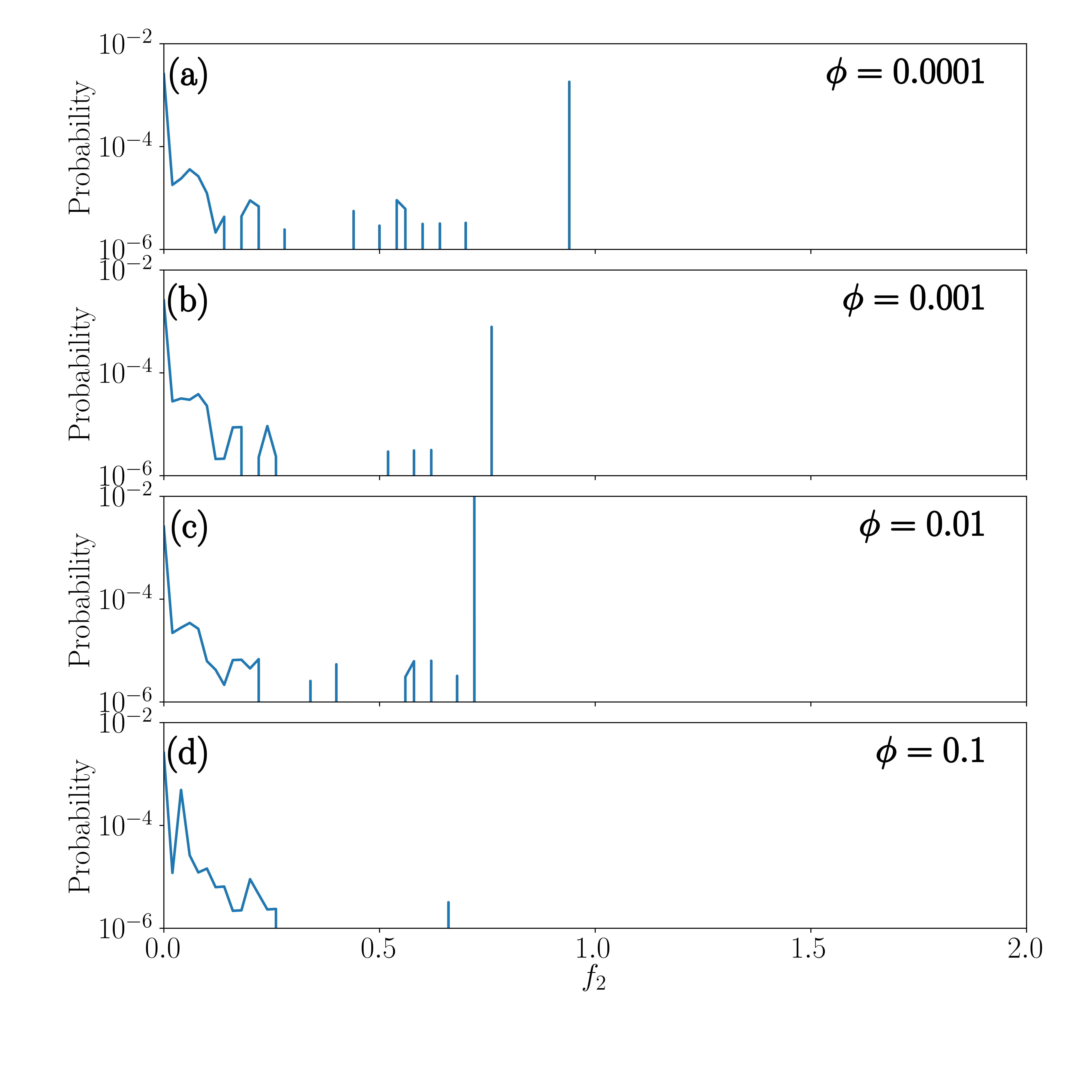}
	\caption{Plots representing the probability of occurrence of extreme events for various values of $\phi$. }
	\label{phicheck}
\end{figure}

We found that upon increasing the value of $\phi$, we can suppress the extreme events that were found beyond $f_2=0.5$. This is evident from Fig.~\ref{phicheck} drawn on the probability of occurrence of extreme events. Such suppression of extreme events can be observed only in the case of non-polynomial mechanical system. In the case of Li\'enard system, extreme events were neither enhanced nor suppressed by the introduction of this new phase.

\subsection{Combined effect of weak second periodic forcing and constant bias}
\label{dfcons}

\begin{figure}[!ht]
	\centering
	\includegraphics[width=0.75\textwidth]{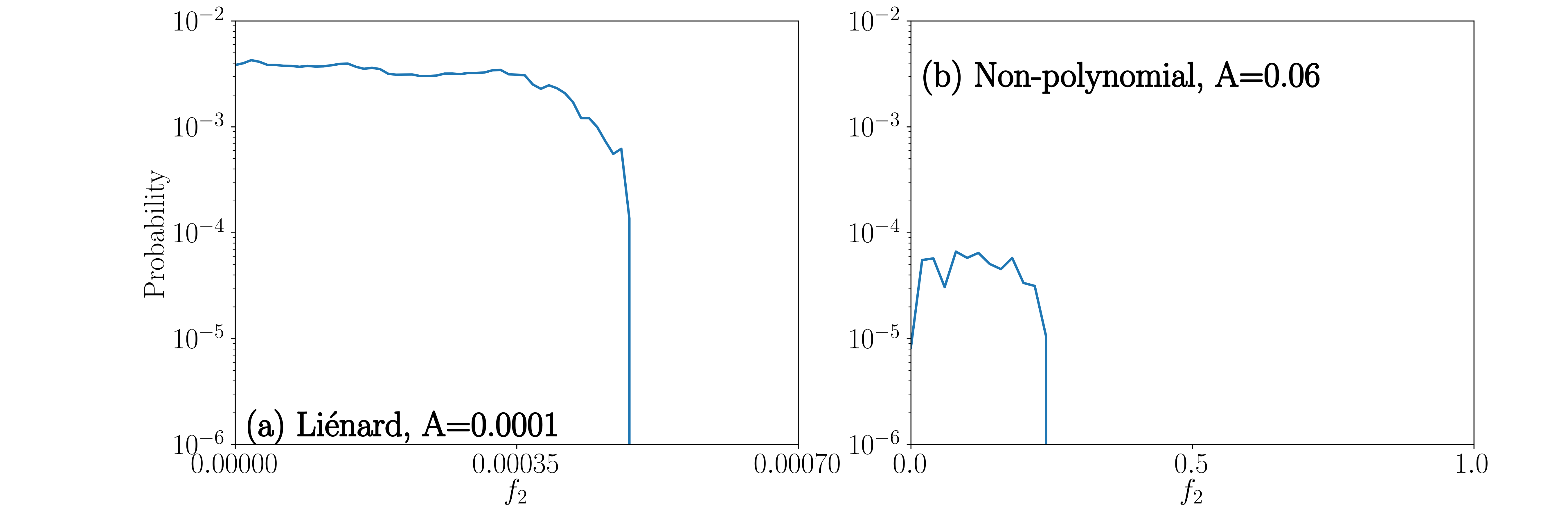}
	\caption{Plots representing the probability of occurrence of extreme events. (a) for Li\'enard system with $A=0.0001$ and (b) for non-polynomial system with $A=0.06$. }
	\label{dual}
\end{figure}

In Fig.~\ref{dual}, we plot the probability of occurrence of extreme events by varying the second forcing strength $f_2$. The value of the constant bias is fixed as $A=0.0001$ and $A=0.06$ for the Li\'enard system and non-polynomial system respectively.

We confirm from Fig.~\ref{dual}, that extreme events gets suppressed even under the combined influence of second periodic forcing and constant bias. In the case of Li\'enard system extreme events get suppressed when $f_2=0.0005$ and in the case of non-polynomial system, extreme events get suppressed when $f_2=0.24$. In circumstances wherein both second periodic force and constant bias are present in the considered system, suppression under combined influence is applicable.

\section{Effect of weak second periodic forcing and constant bias on the parametrically driven non-polynomial system}
\label{paramdfcons}
The emergence of extreme events in a parametrically driven non-polynomial system and its mitigation using single forcing have been studied very recently in \cite{sudharsan}. Equation (\ref{eqn-m-cons}) can be considered as a parametrically driven system, when the angular velocity of the system is driven in the form $\Omega=\Omega_0(1+\epsilon~\mathrm{cos}~\omega_pt)$. The equation of motion with the parametric drive strength $\epsilon$, first forcing $f_1$, second forcing $f_2$ and constant bias $A$ reads now
\begin{eqnarray}
\dot{x}&=&y,\nonumber \\
\dot{y}&=&\frac{-\alpha y-\lambda x y^2 - \omega_0^2 x + \Omega_0^2[2\epsilon~\mathrm{cos}~\omega_pt + 0.5\epsilon^2(1+\mathrm{cos}~2\omega_pt)]x}{1 + \lambda x^2}  \nonumber \\
&&+~\dfrac{f_1 \cos  \omega_1 t + f_2 \cos  \omega_2 t + A}{1 + \lambda x^2}.
\label{eqn-m-param}
\end{eqnarray} 
Here $\omega_1$, and $\omega_2$ are the frequencies of first and second external periodic forcing respectively. System parameters are fixed same as in Sec.~\ref{cons} and we set the values for the new parameters as $\Omega_0^2=6.7$, $\epsilon=0.081$ and $\omega_p=1.0$.

\begin{figure}[!ht]
	\centering
	\includegraphics[width=0.75\textwidth]{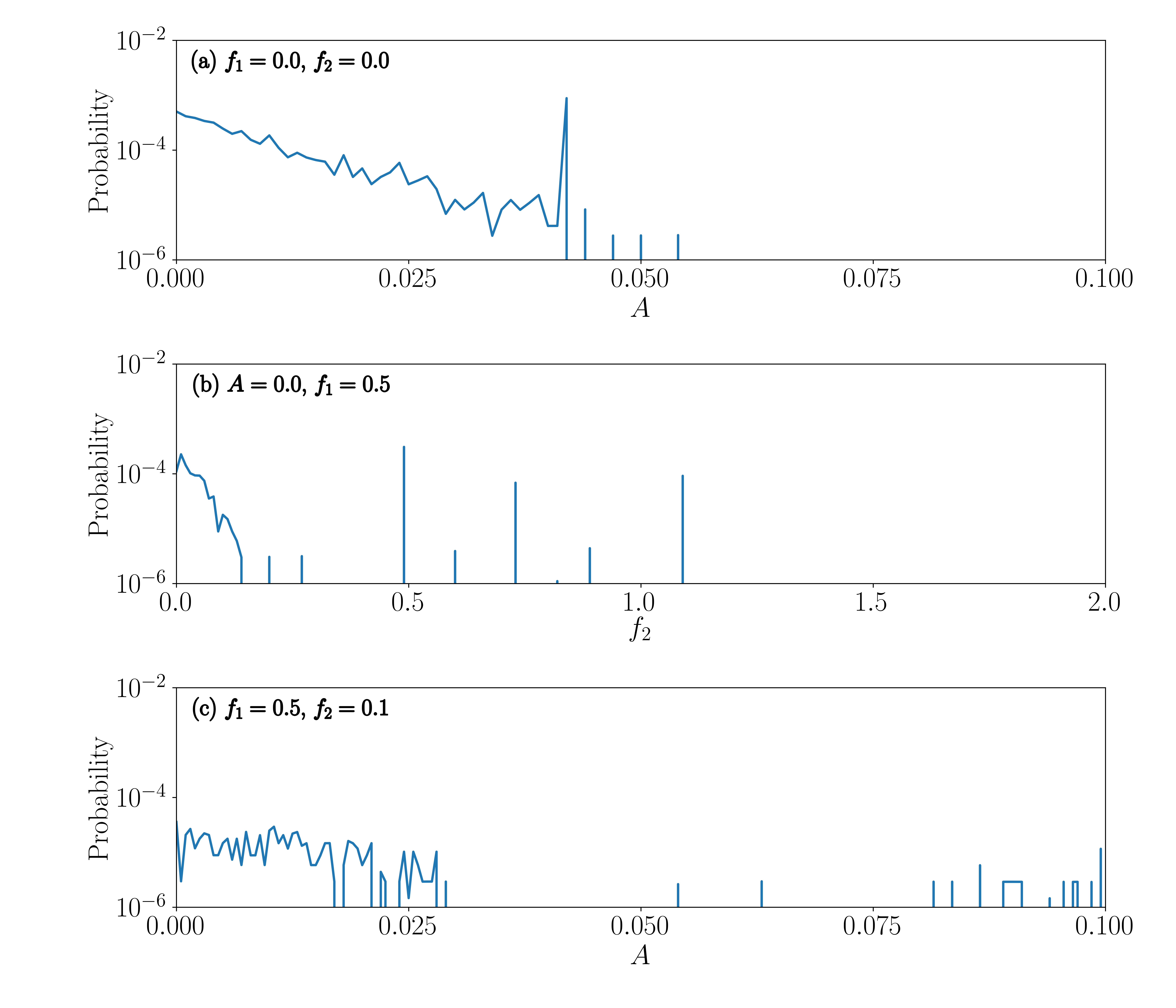}
	\caption{Plots representing the probability of occurrence of extreme events of system (\ref{eqn-m-param}). (a) for $f_1=0.0$, $f_2=0.0$ with varying $A$ (b) for $A=0.0$, $f_1=0.5$ with varying $f_2$, (c) for $f_1=0.5$, $f_2=0.1$ with varying $A$.}
	\label{param}
\end{figure}

Now we check the efficacy of the non-feedback methods on suppressing the extreme events in system (\ref{eqn-m-param}) for three different possibilities, namely (i) parametric drive with constant bias, (ii) parametric drive with dual forcing and (iii) parametric with both dual forcing and constant bias. The results for these three settings are shown in Fig.~\ref{param}.

In all the three cases, in the parametrically driven system (\ref{eqn-m-param}), the suppressing nature of extreme events are explicit. In the system (\ref{eqn-m-param}) also we find that the constant bias suppresses extreme events quicker in comparison with the second periodic forcing. In Fig.~\ref{param}(a), $f_1$ and $f_2$ are fixed to 0 and $A$ is varied. We find that the probability of occurrence of extreme events decreases gradually and becomes zero at $A=0.055$. In the case shown in Fig~\ref{param}(b), $A$ is nullified and $f_1$ is fixed to 0.5 and $f_2$ is varied. We observe that the system decreases in producing extreme events and produces zero extreme events at $f_2=0.15$. Afterwards, extreme events emerge discontinuously, but it is completely suppressed at $f_2=1.1$. Finally, in the case shown in Fig.~\ref{param}(c), $f_1$ and $f_2$ are fixed as 0.5 and 0.1 respectively and $A$ is varied.  In this case also, we notice that extreme events decreases to zero at $A=0.0295$. However, the extreme events emerge later at $A=0.054$ and occur discontinuously after this value. Although the emergence of extreme events occurs later while increasing $A$, the system exhibits zero extreme events for most of the values of $A$. Hence weak second periodic forcing and constant bias suppresses the extreme events in the parametrically driven non-polynomial system also.

\section{Effect of Constant Bias and Second periodic forcing on the multistability.}
\label{sec:multi}

Recently, works have been dedicated to study multistability in chaotic systems \cite{Natiq2019,Rahim2019}. In this direction Li\'enard system (\ref{eqn-l-cons}) has been shown to possess multistablility \cite{Kingston2017,Ouannas2021}. In the following, we investigate the effect of the considered non-feedback methods on multistability. To begin, we consider the Li\'enard system in the following form
\begin{eqnarray}
\dot{x}= y \quad
\dot{y}= -\alpha x y - \gamma x -\beta x^3 + f_1 \mathrm{sin}(\omega_1 t) + f_2\mathrm{sin}(\omega_2 t) + A.
\label{eqn-multi}
\end{eqnarray}
\par For this study, we fix the parameter values as $\alpha=0.0135$, $\beta=0.8111$, $\gamma=-2.65$, $f_1=2.0$. Initially $f_2,~\omega_2$ and $A$ are set to zero while $\omega_1$ is varied in the analysis. The form of the equation and the value of parameters are fixed as in Ref.~\cite{Kingston2017}. 

\begin{figure}
	\centering
	\begin{subfigure}[b]{0.47\textwidth}
		\centering
		\includegraphics[width=\textwidth]{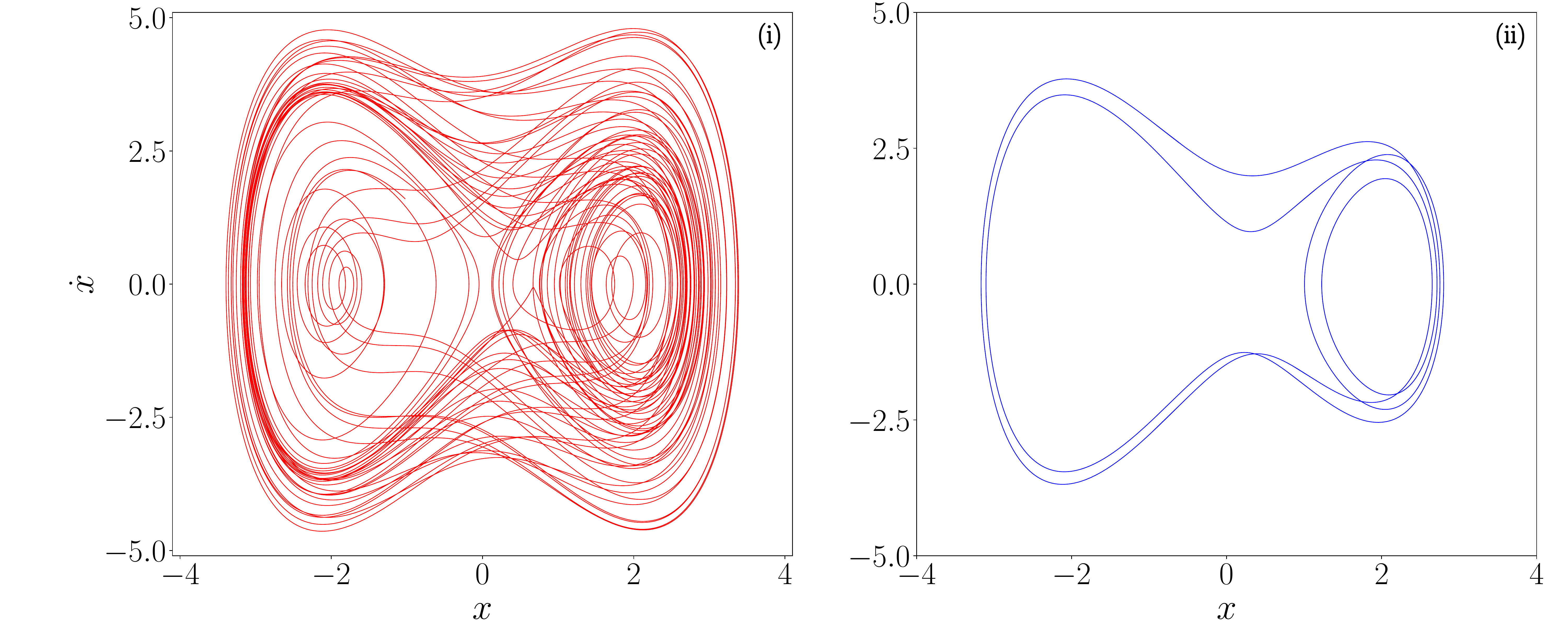}
		\caption{For $\alpha=0.0135$, $\beta=0.8111$, $\gamma=-2.65$, $f_1=2.0$ \& $\omega_1=0.762$. For the set of initial conditions (i) $(x_0,y_0)=(1.8,1.6)$ representing the chaotic nature and (ii) $(x_0,y_0)=(1.67,-2.02)$}
		\label{icphase}
	\end{subfigure}
	\hfill
	\begin{subfigure}[b]{0.47\textwidth}
		\centering
		\includegraphics[width=\textwidth]{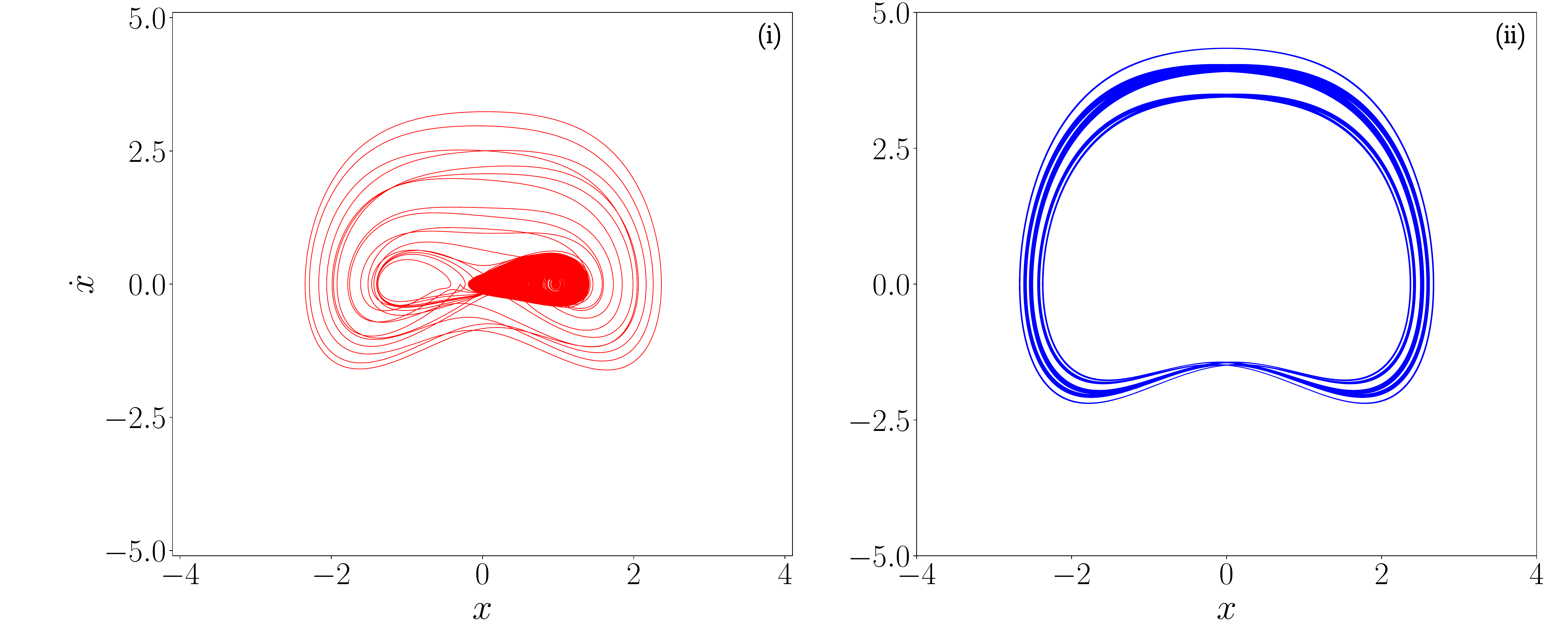}
		\caption{For $\alpha=0.45$, $\gamma=-0.5$, $\beta=0.5$, $f_1=0.2$ $\omega_1=0.7315$ for the set of initial conditions (i) $(x_0,y_0)=(0.5,0.5)$ representing the chaotic nature and (ii) $(x_0,y_0)=(1.67,-2.02)$ representing the quasiperiodic nature }
		\label{icphase1}
	\end{subfigure}
	\caption{Plots representing the phase portraits of system (\ref{eqn-multi}) without constant bias and second periodic forcing. }
	\label{icphasefull}
\end{figure}

Initially, we fix $\omega_1=0.762$ and choose the initial conditions as $(x_0,y_0)=(1.8,1.6)~ \text{and}~ (1.67,-2.02)$ for which the sytem exhibits chaotic and periodic behaviour respectively. The chaotic and the periodic nature of the trajectories for these two initial conditions are brought out through the phase portraits and shown in Fig.~\ref{icphase}. 

\begin{figure}
	\centering
	\begin{subfigure}[b]{0.47\textwidth}
		\centering
		\includegraphics[width=\textwidth]{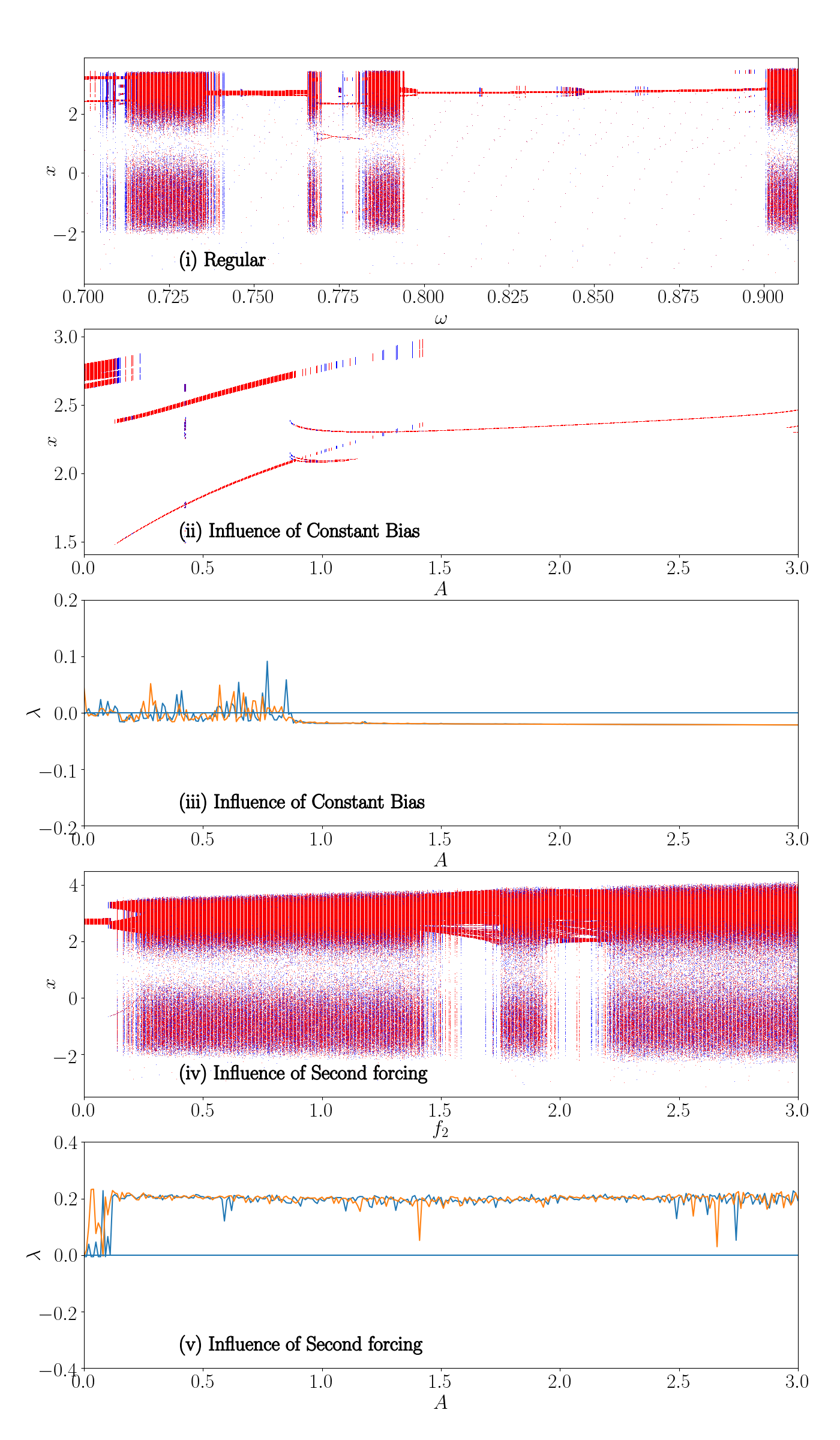}
		\caption{For $\alpha=0.0135$, $\beta=0.8111$, $\gamma=-2.65$, $f_1=2.0$ \& $\omega_1=0.762$.}
		\label{multi}
	\end{subfigure}
	\hfill
	\begin{subfigure}[b]{0.47\textwidth}
		\centering
		\includegraphics[width=\textwidth]{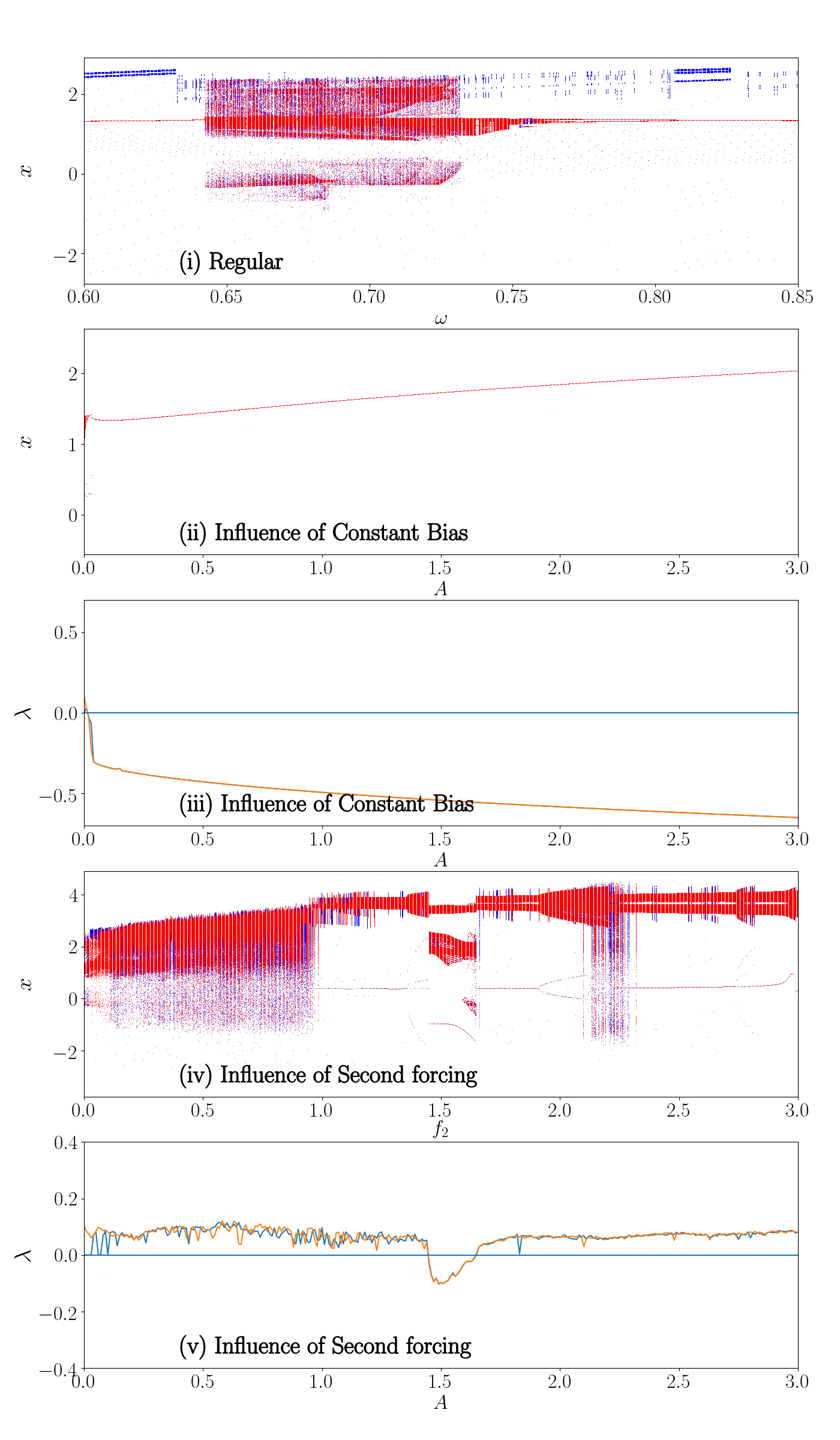}
		\caption{For $\alpha=0.45$, $\gamma=-0.5$, $\beta=0.5$, $f_1=0.2$ $\omega_1=0.7315$}
		\label{multi1}
	\end{subfigure}
	\caption{Plots representing the bifurcation and maximal Lyapunov exponents of system (\ref{eqn-multi}) (i) without the influence of constant bias and second periodic forcing, (ii) \& (iii) under the influence of constant bias and (iv) \& (v) under the influence of the second periodic forcing.}
	\label{multifull}
\end{figure}

In order to analyse the influence of the two non-feedback methods on multistability, we plot in Fig.~\ref{multi} the one parameter bifurcation diagram of the system (\ref{eqn-multi}) by collecting all the maxima. In this figure, blue colour corresponds to the initial conditions $(1.6,1.8)$ whereas red colour corresponds to the initial conditions $(1.67,-2.02)$. In particular, Fig. \ref{multi}(i) is drawn when the system is not influenced either by the constant bias or second periodic forcing. We plot bifurcation diagram and maximal Lyapunov exponents, respectively in Figs. \ref{multi}(ii) \& (iii) when the system is under the influence of constant bias alone and the same in Figs.~\ref{multi}(iv) \& (v) when the system is under the influence of the second periodic forcing alone. From Figs.~\ref{multi}(ii) \& (iii) we can see that under the influence of constant bias, the multistability nature of the system decreases and only periodic orbits prevail for the higher values of $A$. On analysing Figs.~\ref{multi}(iv) \& (v), we observe that  chaotic nature prevails but the multistability nature supresses when we increase the strength of the second forcing $f_2$. From this we conclude that the multistability is brought under control when the system is influenced by the constant bias as well as the second periodic forcing. It is important to note that for the considered choice of parameter the sytem possess a double well potential \cite{Kingston2017}. Now when the system is influenced by a constant bias the system's dynamics is forced to get confined in a single well. This is the reason behind the suppression of chaos as well.    

Now we make the same analysis for the parameter values $\alpha=0.45$, $\gamma=-0.5$, $\beta=0.5$, $f_1=0.7315$. It is for these parameter values that extreme events occur in the Li\'enard system (\ref{eqn-multi}) with neither second periodic forcing nor constant bias. Now we choose the initial conditions as (0.5, 0.5) and (1.67, -2.02) for which the system exhibits chaotic and quasiperiodic behaviour respectively. The phase portraits for these two sets of initial conditions are shown in Fig.~\ref{icphase1}.

The corresponding bifurcation diagram and the maximal Lyapunov exponents for these two sets of initial conditions for the regular, constant bias and second periodic forcing are shown in Fig.~\ref{multi1}. In this case also, the outcome subtantiates the fact that constant bias as well as the double forcing suppresses the multistability nature of the system. The chaotic nature prevails under the influence of the second periodic forcing here as well. 

\begin{figure}[!ht]
	\centering
	\includegraphics[width=1.0\textwidth]{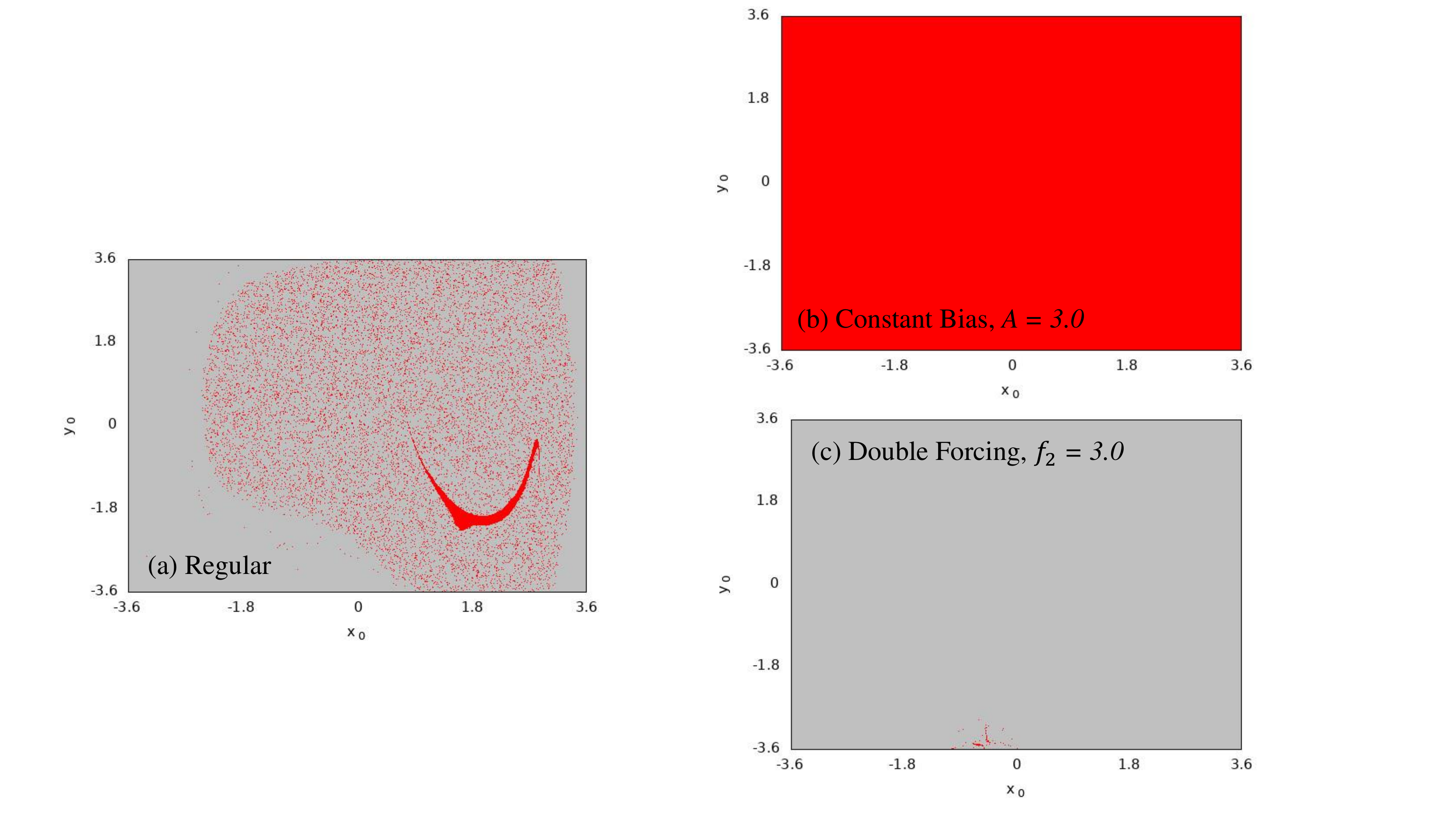}
	\caption{Plots representing the bifurcation of system (\ref{eqn-multi}) with $\omega_1=0.758$. (a) without the influence of constant bias and second periodic forcing, (b) under the influence of constant bias and (c) under the influence of the second periodic forcing.}
	\label{basin}
\end{figure}

In order to substantiate our result further, we plot the basin of attraction in Fig.~\ref{basin} for the same set of parameters as in Fig.~\ref{multi} but with $\omega=0.758$. The value of parameters were chosen as in Ref.~\cite{Kingston2017}. In the basin of attraction grey represents the chaotic regions and red represents the periodic orbits. Initially without constant bias or second periodic forcing there are mixtures of these two states. As we increase the constant bias strength, the periodic nature of the system increases and becomes completely periodic (see Fig.~\ref{basin}(b)) while the nature of the system becomes completely chaotic in the case of increasing second periodic forcing (see Fig.~\ref{basin}(c)). But there are very few regions in the basin of attraction where periodic attractor still prevails which can be confirmed from Fig.~\ref{basin}(c). But this region is very small and hence we can confirm that multistability can be reduced to a considerable extent by the second periodic forcing.

\section{Discussion}
\label{discussion}
The investigations carried out in this paper confirm that constant bias and weak second periodic forcing both can be used as a potentially viable non-feedback method to suppress extreme events. Eventhough we have demonstrated the applicability of suppression of extreme events by considering two well known nonlinear dynamical systems, namely Li\'enard and a non-polynomial mechanical system, the underlying ideas can be extended to a class of nonlinear systems which exhibit extreme events. In the case of Li\'enard system, under the influence of constant bias, extreme events are mitigated by the decrease in large amplitude oscillations. On the otherhand when it is influenced by a weak second periodic forcing, extreme events are mitigated due to the increase in chaotic nature of the system. The increasing chaotic nature is attributed to the increase in large amplitude oscillations in the Li\'enard system. Increase in the large amplitude oscillations correspondingly increases the threshold value thereby making no peaks to cross the threshold. Hence we come across zero extreme events. In the case of non-polynomial system also extreme events get suppressed by the destruction of chaotic attractor under the influence of constant bias. The reason for the suppression of extreme events in non-polynomial system under the influence of second periodic forcing is due to the increase in large amplitude oscillations. We have also verified the effect of constant bias and weak second periodic forcing on the non-polynomial system with parametric drive since very recently the latter system has been shown to exhibit extreme events. We have also confirmed that extreme events get suppressed in this system as well. In the case of constant bias the extreme events gets suppressed for a very small value of the bifurcation parameter value  and in the case of second periodic forcing it happens for a comparatively higher value. 

Upon comparing the suppression of extreme events due to constant bias and weak second forcing, we conclude that the suppression due to constant bias is more effective from the point of view of dynamics. The major difference lies in the fact that in the case of constant bias along with suppression of extreme events, either large amplitude oscillations gets suppressed completely or chaotic attractor is destroyed, whereas in the case of second forcing, large amplitude oscillations persists even after the suppression of extreme events. Presence of large amplitude oscillations are prone to produce extreme events, if the system is influenced by some other external factors such as noise. But if the oscillations are of small amplitude, influence of external factors will not give rise to extreme events. So in this context, suppression of extreme events by constant bias is more effective when compared to the second forcing.

Further it is clear from our analysis that both the constant bias and second periodic forcing supressess the multistability nature that is present in the Li\'enard system. The only difference which we came across is that the constant bias supresses the multistability nature completely whereas the second periodic forcing suppresses the multistability only to a considerable extent.

From the applications point of view, one may consider weak second periodic force as the weak sinusoidal AC current and the constant bias as the fixed DC current.  Eventhough several works have been devoted on the control of extreme events through several feedback methods, the results which we have arrived through this work turns out to be very important from the point of view of reliability. Depending on the requirement, we can choose one of these two methods. Say for example, if we need to eliminate extreme events along with the elimination of chaotic nature and multistability we can opt for the constant bias, whereas if we need to get rid of extreme events without the suppression of chaotic attractor and/or multistability, we can opt for the second periodic forcing. From the implementation point of view it is well known that both the weak AC and  DC current feedbacks are easily implementable in physical systems and in particular when it comes to electronic circuits. Not only implementation but monitoring also becomes very easy and no online monitoring is necessary when we apply such AC and DC currents to systems. So our results turns out to be easily implementable for the control of extreme events in physical, mechanical and electronic systems. This is just a starting point towards the creation of easy control measures for the mitigation of extreme events. The application of weak second periodic forcing and constant bias for the control of extreme events can further be extended to laser based optical systems, climatic systems and oceanic systems. 

\section{Conclusion}
\label{conclusion}
In this work, we have shown that two non-feedback methods, namely constant bias and weak second periodic forcing are potentially viable tools for the mitigation of extreme events in Li\'enard system and in a damped and driven non-polynomial mechanical system. We have demonstrated the suppression of extreme events using the peak PDF which displayed a transition from long tailed distribution (when extreme events occur) to Gaussian like distribution (when no extreme events occurred). We have also verified our results using the probability plot which displayed a gradual decrease to zero in the probability of extreme events. Finally, the $d_{max}$ was calculated and plotted against the bifurcation parameter. Whenever extreme events occur, the value of $d_{max}$ crosses the value $n$ present in the qualifier threshold. Also the introduction of phase $\phi$ to the second periodic force further suppresses the extreme events in the non-polynomial system. Further, we have also considered the combined effect of second periodic forcing and constant bias and found that on fixing the constant bias $A$ and varying $f_2$, extreme events get suppressed. We have consolidated the above results and presented them in a two parameter plot in the $(f_1-A)$ plane and in the $(f_1-f_2)$ plane. Our investigations also reveal that these two non-feedback methods suppresses the extreme events in an effective manner when the angular frequency of the non-polynomial system is found to be parametrically driven. We have also studied the suppression of multistability in Li\'enard system under the action of constant bias and second periodic forcing.

\section*{Acknowledgments}
SS  thanks  the  Department  of  Science  and  Technology (DST), Government of India, for support through INSPIRE Fellowship (IF170319). The work of AV forms a part of a research project sponsored by DST under the Grant No. EMR/2017/002813.  The work of MS forms a part of a research project sponsored by Council of Scientific and Industrial Research (CSIR) under the Grant No. 03(1397)/17/EMR-II. MS also acknowledges the Department of Science and Technology (DST) under PURSE Phase-II for providing financial support in procuring high performance desktop which highly assisted this work.  

\section*{Author Contribution Statement:}
All the authors contributed equally to the preparation of this manuscript. 

\section*{Data Availability Statement }
All data generated or analysed during this study are included in this published article.

\end{document}